\def\mode{1}
\if0\mode
\documentclass[journal,twoside,web]{ieeecolor}
\usepackage{generic}
\else
\documentclass[journal,twoside]{IEEEtran}
\fi
\usepackage[english]{babel}
\usepackage{microtype}
\usepackage{setspace}
\usepackage{blindtext}
\usepackage{siunitx}
\usepackage{titlecaps}
\Addlcwords{or with if of and an to -off off for in into the versus vs subject} % excluded words
\usepackage{comment}

\usepackage[caption=false]{subfig}
\usepackage{graphicx}
\usepackage{fancyhdr} 
\usepackage{color}
\usepackage{epsfig}
\graphicspath{{Images/}}
\usepackage{tikz}
\usetikzlibrary{patterns,fit,matrix}
\usepackage{tikzscale}
\usepackage{scalerel}
\usetikzlibrary{arrows}
\usetikzlibrary{plotmarks}
\usetikzlibrary{svg.path}
\usetikzlibrary{shapes.multipart}
\usepackage{pgfplots}
\usepgfplotslibrary{fillbetween}
\pgfplotsset{compat=newest}
\pgfplotsset{every axis/.append style={
		label style={font=\Large},
		tick label style={font=\large}  
}}
\tikzstyle{int}=[draw, fill=black!10, minimum size=5em,thick]
\tikzstyle{init} = [pin edge={to-,thick,black}]

\usepackage{array}
\usepackage{stfloats}

\usepackage{amsmath,amssymb,amsfonts,amsthm}
\numberwithin{equation}{section}
\usepackage{mathdots}
\usepackage{bbm,dsfont}
\usepackage{relsize}
\usepackage{nicefrac}
\usepackage{bbm}
\usepackage{mathtools}
\usepackage[mathscr]{eucal}

\usepackage[ruled]{algorithm}
\usepackage{algpseudocode}
\usepackage{multirow}
\usepackage{makecell}
\makeatletter
\gdef\Shortstack{\@ifnextchar[\@Shortstack{\@Shortstack[c]}}
\gdef\@Shortstack[#1]#2{%
	\leavevmode
	\vbox\bgroup
	\baselineskip-\p@\lineskip 3\p@
	\let\mb@l\hss\let\mb@r\hss
	\expandafter\let\csname mb@#1\endcsname\relax
	\let\\\@stackcr\setlength{\baselineskip}{#2}%
	\@ishortstack}
\makeatother
\usepackage{booktabs}
\usepackage{tabularx}

\usepackage{enumitem}

\usepackage{float}
\usepackage[absolute]{textpos}

\usepackage{url}
\usepackage{hyperref}
\makeatletter
\let\NAT@parse\undefined
\makeatother
\usepackage[capitalize]{cleveref}
\usepackage{cite}

\newcommand\orcidicon[1]{\href{https://orcid.org/#1}{\includegraphics[scale=0.04]{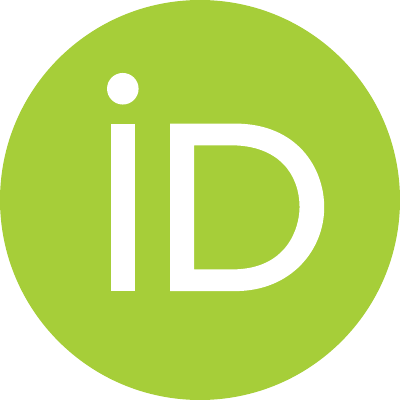}}}

\if1\mode
\definecolor{nblue}{gray}{0}
\definecolor{sectioncolor}{gray}{0}
\definecolor{subsectioncolor}{gray}{0}
\fi

\if0\mode
\newcommand{\startSentence}[1]{{\color{nblue}\sf\itshape\small\bfseries #1}}
\newcommand{\myParagraph}[1]{{\color{nblue}\sf\itshape\small\bfseries {#1.}}}
\else
\newcommand{\startSentence}[1]{{\bfseries #1}}
\newcommand{\myParagraph}[1]{{\bfseries {#1.}}}
\fi

\newcommand{\tradeoff}{centralized-decentralized trade-off\xspace}

\newcommand{\Real}[1]{ { {\mathbb R}^{#1} } }
\newcommand{\Realp}[1]{ { {\mathbb R}^{#1}_+ } }

\DeclareMathOperator*{\argmin}{argmin}

\newcommand{\DeclareAutoPairedDelimiter}[3]{%
	\expandafter\DeclarePairedDelimiter\csname Auto\string#1\endcsname{#2}{#3}%
	\begingroup\edef\x{\endgroup
		\noexpand\DeclareRobustCommand{\noexpand#1}{%
			\expandafter\noexpand\csname Auto\string#1\endcsname*}}%
	\x}

\DeclareAutoPairedDelimiter{\ceil}{\lceil}{\rceil}
\DeclareAutoPairedDelimiter{\floor}{\lfloor}{\rfloor}
\newcommand{\tr}[1]{\mbox{tr}\left(#1\right)}
\newcommand{\gauss}{\mathcal{N}}

\renewcommand{\mod}{\mathrm{ mod}}
\newcommand{\consMatrix}{\dfrac{\mathds{1}_N\mathds{1}_N^{\top}}{N}}
\newcommand{\e}{\mathrm{e}}

\newcommand{\x}[2]{x_{#1}(#2)}

\newcommand{\xtilde}[2]{\tilde{x}_{#1}(#2)}

\newcommand{\xbar}[2]{\bar{x}_{#1}(#2)}

\newcommand{\ztilde}[2]{\tilde{z}_{#1}(#2)}
\newcommand{\zbar}[2]{\bar{z}_{#1}(#2)}

\newcommand{\noise}[2]{w_{#1}(#2)}
\newcommand{\noisebar}[2]{\bar{w}_{#1}(#2)}
\newcommand{\noisetilde}[2]{\tilde{w}_{#1}(#2)}

\renewcommand{\u}[2]{u_{#1}(#2)}

\newcommand{\meas}[3]{y_{#1,#2}(#3)}
\newcommand{\delayn}{\tau}
\newcommand{\taun}{\tau_n}

\newcommand{\gvel}{\eta}
\newcommand{\gpos}{\lambda}
\newcommand{\opteig}{\lambda^*}
\newcommand{\setstable}{\mathcal{S}}

\newcommand{\varx}[2]{\sigma^{2}_{\textit{#2}}\left(#1\right)}
\newcommand{\var}{\sigma^{2}}

\if0\mode
\newtheoremstyle{colorthm}{\topsep}{\topsep}{\normalfont}{1em}{\color{nblue}\sf\itshape\small\bfseries}{:}{ }{\thmname{#1}\thmnumber{ #2}{\thmnote{ (#3)}}}
\theoremstyle{colorthm}
\else
\theoremstyle{plain}
\fi
\newtheorem{thm}{Theorem}
\newtheorem{cor}{Corollary}
\newtheorem{prop}{Proposition}

\newtheorem{prob}{Problem}
\newtheorem{ass}{Assumption}
\newtheorem{rem}{Remark}

\addto\captionsenglish{}
\addto\captionsenglish{}

\newcommand{\canOmit}[1]{{\color{gray}#1}\xspace}
\newcommand{\review}[1]{{\color{black}#1}}
\newcommand{\revision}[1]{{\color{black}#1}}
\newcommand{\revisiontwo}[1]{{\color{black}#1}}
\newcommand{\red}[1]{\textcolor{red}{#1}}
\newcommand{\done}[1]{{\renewcommand\tcb[1]{{##1}}#1}}
\newcommand{\tcb}[1]{{\color{blue} #1}} % Mihailo

\newcommand{\blue}[1]{{\color{black} #1}}
\newcommand{\mjmargin}[1]{\marginpar{\color{red}\tiny\ttfamily{MJ:} #1}}

\newcommand{\linkToPdf}[1]{\href{#1}{\blue{(pdf)}}}
\newcommand{\linkToPpt}[1]{\href{#1}{\blue{(ppt)}}}
\newcommand{\linkToCode}[1]{\href{#1}{\blue{(code)}}}
\newcommand{\linkToWeb}[1]{\href{#1}{\blue{(web)}}}
\newcommand{\linkToVideo}[1]{\href{#1}{\blue{(video)}}}
\newcommand{\linkToMedia}[1]{\href{#1}{\blue{(media)}}}
\newcommand{\award}[1]{\xspace} % omit awards
\newcommand{\eg}{\emph{e.g.,}\xspace}
\newcommand{\ie}{\emph{i.e.,}\xspace}
\addto\extrasenglish{}

\addto\extrasenglish{}
\addto\extrasenglish{}
\addto\extrasenglish{}

\title{{\titlecap{can decentralized control outperform centralized?} \\ \titlecap{the role of communication latency}}}

\author{Luca~Ballotta\textsuperscript{\orcidicon{0000-0002-6521-7142}}, %~\IEEEmembership{Graduate~Student~Member,~IEEE},\\%
	Mihailo~R.~Jovanovi\'c\textsuperscript{\orcidicon{0000-0002-4181-2924}},~\IEEEmembership{Fellow,~IEEE}, and %
	Luca~Schenato\textsuperscript{\orcidicon{0000-0003-2544-2553}},~\IEEEmembership{Fellow,~IEEE}%
	\thanks{
		This work has been partially supported
		by the Italian Ministry of Education, University and Research (MIUR) through
		the PRIN project no. 2017NS9FEY entitled ``Realtime Control of 5G Wireless Networks'', and through
		the initiative "Departments of Excellence" (Law 232/2016), and by
		the US National Science Foundation (NSF) under Awards ECCS-1708906 and ECCS-1809833.
		Views and opinions expressed in this work are of the authors and may not reflect those of the funding institutions.}%
	\thanks{Luca Ballotta and Luca Schenato are with the Department of Information Engineering, University of Padova, 35131 Padova, Italy
		(e-mail: ballotta@dei.unipd.it; schenato@dei.unipd.it)}%
	\thanks{Mihailo R.\ Jovanovi\'c is with the Ming Hsieh Department of Electrical and Computer Engineering,
		University of Southern California, Los Angeles, CA 90089 USA (e-mail: mihailo@usc.edu)}
}

\if0\mode
\fi

\begin{document}
	
	\if1\mode
	\begin{textblock}{20}(-2,0.05)
		\footnotesize
		\centering
		\setstretch{1}
		This article has been accepted for publication on the IEEE Transactions on Control of Network Systems.\\
		Please cite the paper as: L. Ballotta, M. R. Jovanovi\'c, and L. Schenato,\\
		“{\titlecap{can decentralized control outperform centralized?} \titlecap{the role of communication latency}}”,\\
		IEEE Transactions on Control of Network Systems, 2023.\\
		%		Link to article abstract: \url{}
	\end{textblock}
	\fi
	
	\if0\mode
	\bstctlcite{MyBSTcontrol}
	\fi
	
	\maketitle
	%!TEX ROOT = ../centralized_vs_distributed.tex

\begin{abstract}
	\done{
	In this paper,
	we \tcb{examine the influence} of 
	communication latency on
	\tcb{performance} of networked control systems.
	\tcb{Even though} \review{distributed} \tcb{control} architectures \tcb{offer
		advantages in terms of communication}, maintenance costs, \tcb{and} scalability, 
	\tcb{it is an open question how communication latency
		that varies with network topology
		influences closed-loop performance}.
	\tcb{For networks in which delays increase with the number of links, we
		establish the existence of a fundamental performance trade-off that arises from 
		control architecture.
		In particular, we utilize consensus dynamics with single- and 
		double-integrator agents to show that,
		if delays increase fast enough, 
		a sparse controller with nearest neighbor interactions can outperform 
		the centralized one with all-to-all communication topology.} 
	}
	
\end{abstract}
	%!TEX ROOT = ../../centralized_vs_distributed.tex

\section{Introduction}\label{sec:intro}

%% We all agree that decentralized control is the way to go
\IEEEPARstart{I}{t is} widely accepted that
modern multi-agent systems cannot rely on centralized control architectures.
This conclusion stems from issues related to gathering
all decision making to a central node,
ranging from lack of robustness and failures proneness,
to maintenance costs,
and communication overhead.
%While such issues may be neglected for small systems,
%they represent critical bottlenecks for
Indeed, large-scale networks
have experienced a net shift
towards decentralized and distributed architectures~\cite{JOVANOVIC201676,8340193}.
Moreover, the recent deployment of powerful communication protocols for massive networks,
\eg 5G~\cite{BIRAL20151,li20185g},
and advances in embedded electronics~\cite{8600375,nvidia},
as well as in algorithms for low-power devices (\eg TinyML~\cite{warden2019tinyml}),
which allow to spread computational tasks across network nodes
according to edge- and fog-computing paradigms~\cite{shi2020joint,shi2016edge,yi2015survey},
are making such networked systems grow at unprecedented scale,
further stressing the importance of distributed controller architectures.

%% Latency sucks but we cannot help it much
A challenging issue in large-scale wireless network systems is
the latency arising from channel constraints, such as
limited bandwidth or packet retransmissions.
%, on the one hand,
%and from limited computational resources at the agents, on the other hand.
To address this problem,
research efforts have been moving towards two main directions.
					%!TEX ROOT = ../../centralized_vs_distributed.tex

%% Control literature
%Many works in the literature support decentralized architectures with various arguments.
\startSentence{Related work in control theory} deals with control design for
distributed architectures,
where classical methods,
such as LQG or $ \mathcal{H}_2 $/$ \mathcal{H}_{\infty} $ control,
\review{require an all-to-all information exchange which is infeasible for large-scale systems.}

%the design % in the presence of latency % in the presence of delays
A \review{large body of work} focuses on stability,
%and the identification of convex structured problems.
\eg~\cite{ren2017finite,SUN2021419} are concerned with finite-time delay-dependent stability of
discrete-time systems,
\cite{BEREZANSKY2015605} finds sufficient conditions for uniform stability
of linear delay systems,
\revision{\cite{MUNZ20101252} characterizes stability and consensus conditions with homogeneous and heterogeneous feedback delays},
and~\cite{Chehardoli2019,8844785} analyze consensus and error compensation for vehicular platoons.
Another line of work deals with maximizing performance for structured controllers,
\revision{\eg~\cite{8358743,Michiels2016,8430769} study $ \mathcal{H}_2 $-norm minimization for time-delay network systems,}
\cite{SOUDBAKHSH2017171} proposes a cyber-physical architecture with LQR for wide-area power systems,
\cite{MORATO202178} develops a procedure for time-varying dead-time compensation
by adapting the Filtered Smith Predictor,
and~\cite{9137405} investigates sensor-and-processing selection for optimal estimation in star networks.

A more recent trend is optimizing the controller architecture.
%namely, the communication network.
For large-scale systems,
this means sparsifying the structure
to enhance communication and scalability.
This is achieved by introducing penalty terms %in the cost function
to trade performance for controller complexity~\cite{6497509,dorjovchebulTPS14,7835692,9216852,7347386,BAHAVARNIA201710395,7378905,ANDERSON2019364}.
In particular,~\cite{7378905} proposes the \textit{Regularization for Design},
addressing optimization of communication links,
\revision{while~\cite{ANDERSON2019364} investigates communication locality
and its relation to control design within the \textit{System Level Synthesis}.}

%% Optimization literature
\startSentence{Related work in optimization theory} is concerned
with minimization of distributed cost functions,
which are only partially accessible at each agent.
%and possibly reveal themselves overtime.
A large body of literature has been devoted to study suitable algorithms,
a short list of which is represented by
\review{\cite{8027140,XIAO200733,8015179,7807315,7472453,mogjovTCNS18}}.
In particular, a line of work has been concerned specifically with the design
of algorithms in the presence of communication delays,
the main issues being related to convergence conditions.
For example,~\cite{6120272,6571230,7994706,ZONG2019412,garcia2016periodic}
study consensus of multi-agent systems with additive or multiplicative time-delays
under various network topologies and agent dynamics.
%showing how its convergence depends on such model parameters.
This approach usually the communication network be given
and focuses on the information exchange and processing by the agents
from an optimization standpoint.
					%!TEX ROOT = ../../centralized_vs_distributed.tex

\myParagraph{Addressed Problem} %{Despite such a bulky literature},
\done{\tcb{Even though} both control design for delay-dependent dynamics 
and \tcb{design of controller architectures are well-studied topics, it remains unclear how} {\em network connectivity affects \tcb{the closed-loop performance in the presence of architecture-dependent communication latency.}} 
\tcb{When the total available bandwidth does not increase with the size of the network~\cite{garcia2016periodic}
	or when multi-hop communication is used among low-power devices~\cite{gupta2010delay},}
the number of active communication links may %cause a non-negligible variation of such a latency.
\revision{affect such latency in non-negligible way}.
\tcb{\revision{In this case, it is important that the control design takes into account increase in delays 
%	with the number of links
	when new communication links are introduced.}}}

\done{
Such an approach is conceptually different from
the \tcb{approaches used} in literature.
On one hand,
delay-aware control designs \revision{such as~\cite{MUNZ20101252,8430769}}
\revision{assume a fixed controller architecture
and either target optimization of the feedback gains
or evaluate stability with respect to gains and/or delays}.
On the other hand,
architecture designs such as~\cite{7378905,7347386}
\revision{do not quantify the impact of architecture-dependent delays on performance,
but explicitly force sparsity by
adding a regularization term that penalizes controller complexity to delay-free performance metrics}.
In fact,
while the fully connected architecture is avoided
because of practical limitations,
it is usually regarded as an upper bound for performance~\cite{JOVANOVIC201676}.}
\revision{To the best of our knowledge, the only works where architecture-dependent delays are used to 
	compute the performance metric are~\cite{gupta2010delay,gupta2011delay},
	where the authors study how transmission power affects convergence rate of consensus.}

\revision{
	We study class of static feedback policies in which control action is formed 
	by utilizing delayed measurements from a limited number of nodes within a network. 
	Impact of similar type of controller architectures 
	on mean-square performance of delay-free stochastically forced consensus, 
	synchronization, 
	and vehicular formation networks has been studied in the literature~\cite{bamjovmitpat12,mogjovTCNS18},
	and our objective is to understand influence of delays on performance
	trade-offs induced by such localized controller architectures relative to centralized ones. 
	Identifying similar trade-offs within other classes of localized control policies 
	(including System Level Synthesis) is a relevant open question which is outside the scope of the current study. 
}

\revision{\myParagraph{Original Contribution}
We aim to bridge the two domains of delay-aware control and architecture design by
quantifying how the latter % with the that appear in the dynamics of the controlled system
affects performance under architecture-dependent communication delays. %is affected by variations in such delays.
%\revision{Specifically, we are interested in \textit{the optimal architecture}.
We address two key challenges. 
First,
we focus on \textit{optimal performance},
whereby \emph{stability} is a prerequisite to control design
needed to provide a bounded cost function. % explodes in the absence of stability.
Hence,
we derive stability conditions that are instrumental to an optimal control design problem.
Second,
we aim to identify the \textit{optimal controller architecture} under delays and quantify fundamental performance trade-offs.
Towards this goal, to circumvent the discrete nature of graphs, 
we work our way through two stages:
first, 
we parametrize each architecture with a parameter $ n $
which characterizes both number of links and delay associated with that architecture,
and show how to compute the optimal controller for a given $ n $. 
We then compare the optimal performance obtained for different values of $ n $, 
which allows us to fairly establish which architectures provide the best closed-loop performance. 
In contrast to~\cite{gupta2010delay,gupta2011delay},}
\revision{
we examine mean-square performance of stochastically forced networks, 
study generic delay functions, 
and address optimal design of feedback gains for different controller architectures.}

%indeed, this approach ties the network connectivity
%directly to the system performance,
%this causes densely connected topologies
%to naturally degrade the performance,
%with no need of adding somewhat abstract penalty terms.
%Indeed, such an approach tightly bounds the controller architecture
%with the dynamic-related performance achieved by the controlled system.
					%!TEX ROOT = ../../centralized_vs_distributed.tex

%\subsection{Preview of key results}\label{sec:contribution}

%Our contribution is twofold.
%	\mjmargin{not sure what light blue specifies}
\if0\mode
\newpage
\fi
\myParagraph{Preview of Key Results}
{We utilize undirected graphs with %ring topology with 
single- and double-integrator agent dynamics to examine fundamental performance limitations in networked systems with architecture-dependent communication delays. 
By exploiting convexity of a minimum-variance control design problem with respect to the feedback gains}, %circulant structure of the underlying matrices 
we demonstrate that the choice of controller architecture
has profound impact on network performance in the presence of delays. 
\revision{In particular, when the delays increase fast enough with the number of links,
sparse topologies can outperform highly connected ones.}
%even though they limit global information exchange.
%This allows to keep the optimization relatively simple
%and to hold the focus on the main matter under investigated,
%which in the following we shall call \textit{\tradeoff}.

%As discussed previously,
%our approach makes network-dependent latency
%directly reflect on the cost-to-go.
%%\ie with different network topologies.
%Contrary to the conventional wisdom that
%the fully connected control architecture leads to
%optimal performance,
%we show that, % a trade-off arises
%if the delays increase fast enough with the number of links,
%%\eg as a result of delayed feedback,
%highly connected topologies perform worse than
%sparse architectures.
%%in general.

\begin{figure}
	\centering
	\includegraphics[width=.67\linewidth]{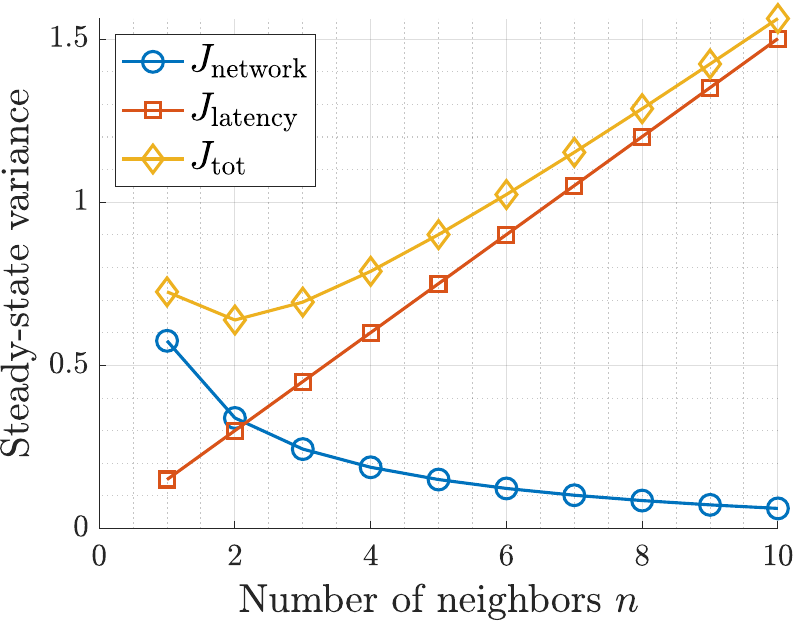}
	\caption{Steady-state variance $ J_{\textrm{tot}}(n) $ versus number of neighbors. %~\eqref{eq:trade-off}.
		The variance is the sum of two costs:
		$ J_{\textrm{network}}(n) $ represents impact of control architecture,
		while
		%		and decreases with denser control architectures.
		$ {J}_{\textrm{latency}}(n) $ is due to the delays affecting the dynamics.
		%		which increases when the addition of links induces longer communication delays.
	}
	\label{fig:trade-off}
\end{figure}

\done{
\tcb{We show that the steady-state variance of a stochastically forced network,
	$ J_{\textrm{tot}}(n) $, can be represented by a sum of two monotone functions of the number of neighbors $ n $ (\autoref{fig:trade-off}),
}}
	\iffalse
	a product of two monotone functions of the number of neighbors $ n $,
\begin{equation}\label{eq:trade-off-mul}
	\tcb{J_{\textrm{tot}}(n) = J_{\textrm{network}}(n)\cdot J_{\textrm{latency}}(n).}
\end{equation}

Here, 
$ J_{\textrm{network}}(n) $ quantifies impact of control architecture and
$ J_{\textrm{latency}}(n) $ determines influence of communication latency on network performance.
While $ J_{\textrm{network}}(n) $ decreases with $ n $ and is minimized by a fully-connected centralized architecture,
$ J_{\textrm{latency}}(n) $ increases with $n$.
This demonstrates the presence of a fundamental trade-off:
on one hand, 
feedback control takes advantage of dense topologies that benefit from information sharing but,
on the other, 
many communication links induce long delays which has negative impact on network performance.

Furthermore, if the delays are sublinear in $ n $,
we show that the trade-off can be also decomposed additive-wise (see ~\autoref{fig:trade-off}):
\fi

\done{\begin{equation}\label{eq:trade-off}
	{\tcb{J_{\text{tot}}(n) = J_{\textrm{network}}(n) + J_{\textrm{latency}}(n)}.}
\end{equation}
Here, 
$ J_{\textrm{network}}(n) $ quantifies impact of control architecture and
$ J_{\textrm{latency}}(n) $ determines influence of communication latency on network performance.
While $ J_{\textrm{network}}(n) $ decreases with $ n $ and is minimized by a fully-connected centralized architecture,
$ J_{\textrm{latency}}(n) $ increases with $n$.
This demonstrates the presence of a fundamental trade-off:
on one hand,
feedback control takes advantage of dense topologies that enhance information sharing but,
on the other hand, many communication links induce long delays which have negative effect on performance.

While~\eqref{eq:trade-off} can be derived analytically
\tcb{for ring topology with continuous-time, single-integrator dynamics},
	our
	\tcb{computational experiments} show that
	a \tcb{similar \textit{\tradeoff} can be observed}
	\revision{\tcb{for general undirected topologies}
	and with double-integrator and discrete-time agent dynamics}. 
	\review{Furthermore,
	in some cases, decentralized architecture with nearest neighbor information exchange
	provides optimal performance.}}

					%!TEX ROOT = ../../centralized_vs_distributed.tex

%\subsection{Paper outline}\label{sec:outline}

\begin{table}
	\centering
	\caption{Theoretical tools (italic) and technical results (roman).}
	\label{tab:results}
	\footnotesize
	\begin{tabular}{|c|c|c|c|}
		\hline
		 & \textbf{Model} & \textbf{Stability} & \textbf{Variance} \\
		\hline
		\multirow{5}{*}{\shortstack{\textbf{Cont.} \\\textbf{time} \\\textbf{(CT)}}} & 
										\makecell{Single int. \\ \eqref{eq:cont-time-single-int-model},\eqref{eq:prop-control}} & 
										\makecell{\textit{Scalar SDDEs~\cite{KuchlerLangevinEqs}} \\ 
											Closed form~\eqref{eq:cont-time-single-int-variance-condition}} &
										\makecell{\textit{Scalar SDDEs~\cite{KuchlerLangevinEqs}} \\
											Closed form~\eqref{eq:cont-time-single-int-steady-state-variance}} \\
		\cline{2-4}
									& \makecell{Double int.\\ \eqref{eq:cont-time-double-int-model}--\eqref{eq:control-input-PD}} & 
									\makecell{\textit{Exponential} \\ \textit{polynomials~\cite{BAPTISTINI1997259}} \\ 
										\textit{SDDEs~\cite{datko1978procedure,wangBoundedness}} \\
										{Implicit}~\eqref{eq:cont-time-double-int-stability-condition}} & 
									\makecell{\textit{SDDEs~\cite{wangBoundedness}, time-}\\
										\textit{scale separation~\cite{khalil2002nonlinear}}\\
										{Integral form}~\eqref{eq:2nd-order-cont-ss-variance}\\
										Approximated~\eqref{eq:x-dynamics-1st-order-approximation}} \\
		\hline
		\multirow{4}{*}{\shortstack{\textbf{Disc.} \\ \textbf{time} \\\textbf{(DT)}}} &  
										\makecell{Single int. \\ \eqref{eq:disc-time-single-int-model},\eqref{eq:prop-control}} & 
										\makecell{\textit{Root locus~\cite{Westphal2001}}\\ 
											Closed form~\eqref{eq:disc-time-single-int-stability-condition}} &
										\makecell{\textit{Moment matching w/} \\
											\textit{Yule-Walker eqs.~\cite{yuleWalkerEqs}} \\
											{Recursive}~\eqref{eq:disc-time-single-int-moment-matching-eqs},\eqref{eq:disc-time-single-int-variance-explicit}} \\
		\cline{2-4}
									&  \makecell{Double int. \\ \eqref{eq:disc-time-double-int-model}} & 
									\makecell{\textit{Jury criterion~\cite{Jury}} \\ 
										{Closed form}~\eqref{eq:disc-time-double-int-characteristic-polinomial}} &
									\makecell{\textit{Moment matching w/} \\
										\textit{Yule-Walker eqs.~\cite{yuleWalkerEqs}} \\
										{Closed form}~\eqref{eq:disc-time-double-int-moment-matching-eqs}} \\
		\hline
	\end{tabular}
\end{table}

\myParagraph{\titlecap{Paper outline}}
In~\autoref{sec:setup} we describe models for communication and controller architecture
and formulate the minimum-variance control design problem.
\revision{While we first utilize ring topology to provide analytical insight
we also demonstrate that our framework can be extended to general undirected topologies; see~\autoref{sec:generic-topology}.}\linebreak
\revision{In Sections~\ref{sec:cont-time}--\ref{sec:cont-time-single-int-control-design}, we %solve the problem for continuous-time agent dynamics,
lay the ground for our main result.
In~\autoref{sec:cont-time}, 
we derive conditions for mean-square stability 
and compute the steady-state variance of continuous-time stochastically forced systems
using Stochastic Delay Differential Equations (SDDEs). 
In~\autoref{sec:cont-time-single-int-control-design},
we prove that the control design problem is convex} \review{and in 
\revision{\autoref{sec:numerical-results} we present} our main results:
by numerically computing the optimal controller gains,
we show that the closed-loop performance is optimized by sparse architectures.
Furthermore, 
we derive analytical expression~\eqref{eq:trade-off} for continuous-time single-integrator dynamics 
which demonstrates that the minimizer is in general nontrivial.}
To address wireless communication,
we study discrete-time systems in~\autoref{sec:disc-time}
and show that the fundamental behavior of the system does not change.
\autoref{tab:results} summarizes our technical results and the theoretical tools used throughout the paper.
\revision{Apart from classical control techniques such as the Jury stability criterion,
we also leverage more unconventional tools from mathematical literature,
such as exponential polynomials~\cite{BAPTISTINI1997259}.}
Concluding remarks are given in~\autoref{sec:conclusion}.
	%!TEX ROOT = ../centralized_vs_distributed.tex

\section{Problem Setup}\label{sec:setup}

We consider {an undirected network} with $ N $ agents 
{in which the state of the $ i $th agent at time $ t $ is given by $ \xbar{i}{t}\in\Real{} $ with the control input $ \u{i}{t}\in\Real{} $.}
For {notational} convenience,
we {introduce} the aggregate {state of the} system  $ \xbar{}{t} $ and {the}
aggregate control input $ \u{}{t} $ {by stacking states and control inputs of each subsystem} $ \xbar{i}{t} $ and $ \u{i}{t} $, 
respectively.

\iffalse
\begin{rem}[Network topology]
	While in the first part of the paper we focus on circular formations for the sake of analysis and ease of presentation,
	the control design can be readily extended to generic undirected topologies.
	We discuss theoretical guarantees in~\autoref{sec:generic-topology} and observe with computational experiments in~\autoref{sec:numerical-results}
	that the \tradeoff holds regardless of the specific topology. % at hand.
\end{rem}
\fi

\myParagraph{Problem Statement}
The agents aim to reach consensus towards a common state trajectory. 
The $i$th component of the vector $ \x{}{t} \doteq \Omega\xbar{}{t} $ represents 
the mismatch between the state of agent $ i $ and the average network state at time $ t $\revision{\cite{bamjovmitpat12}},
where
\begin{equation}\label{eq:error-matrix}
	\Omega \doteq I_{N}-\consMatrix
\end{equation}
and $ \mathds{1}_N \in\Real{N} $ is the vector of all ones,
such that $ \Omega\mathds{1}_N=0 $.
%\red{The target consensus vector is defined as $ \x{m}{t} \doteq \xbar{}{t} - \x{}{t} $.}

\done{
	\tcb{\myParagraph{Ring Topology}
	We focus on ring topology to obtain analytical insights about 
	optimal control design and fundamental performance trade-offs in the presence of communication delays. 
	While some of our notation is tailored to such topology (\eg see equations~\eqref{eq:meas} and~\eqref{eq:feedback-matrix}), 
	in~\autoref{sec:generic-topology} we discuss extension of the optimal control design to generic undirected networks 
	and complement these developments with computational experiments in~\autoref{sec:numerical-results}.}
}

%\myParagraph{\titlecap{Communication model}}
%The agents communicate through a {shared wireless channel}.
%Data are exchanged through a {shared wireless channel} in a symmetric fashion.
%Agent $ i $ communicates with
%\red{$ n $ pairs of agents,
%both agents in each such pair being at equal distance from $ i $}
%\tcb{its} $ 2n $ closest neighbors \tcb{in ring topology.}
%Also, we make the following assumption
%to address channel constraints.

\begin{ass}[Communication model]\label{ass:hypothesis}
	Data are exchanged through a shared wireless channel in a symmetric fashion.
%	\tcb{its} $ 2n $ closest neighbors \tcb{in ring topology.}
	\revisiontwo{Agent $ i $ receives state measurements from
	all agents within $ n $ communication hops.}
	All measurements are received with delay $ \taun \doteq f(n) $
	where $ f(\cdot) $ is a positive increasing sequence.
	\revisiontwo{In particular,
		in ring topology,
		agent $ i $ receives state measurements from the
		$ 2n $ closest agents,
		that is,
		from the $ n $ pairs of agents at distance $ \ell = 1,\dots,n $,
		with $ 1\le n<\nicefrac{N}{2} $.}\footnote{\revision{
%	where both agents in each such pair are at equal distance $ \ell $ from $ i $. % in the ring topology.
%%	located $ \ell $ positions ahead and behind in the formation,.
%	In what follows,
%	without loss of generality,
%	we assume that such $ n $ agent pairs coincide with the
%	$ 2n $ closest agents in ring topology,
%	and that each pair is at distance $ \ell = 1,\dots,n<\nicefrac{N}{2} $.
	\revisiontwo{For example,
	$ n = 1 $ corresponds to nearest-neighbor interaction in ring topology
	and $ n = \floor{\nicefrac{(N-1)}{2}} $ to all-to-all communication topology.}}}
%	Also, each agent measures its own state instantaneously.
\end{ass}

\revision{\begin{rem}[Architecture parametrization]\label{rem:architecture-param}
	Parameter $ n $ will play a crucial role throughout our discussion. 
	In particular,
	we will use it to (i) evaluate the optimal performance %that can be attained 
	for a given
	budget of links
	\revisiontwo{(see~\cref{prob:variance-minimization})};
	and to (ii) compare optimal performance of different control architectures.
	In the first part of the paper, 
	we examine circular formations and
	$ n $ represents how many neighbor pairs communicate with each agent.
	For \linebreak general undirected networks,
	$ n $ determines the number of communication hops for each agent.
	In general,
	$ n $ characterizes sparsity of a controller architecture:
	sparse controllers correspond to small $n$ while highly connected ones to 
	large $ n $.
\end{rem}}

%\begin{rem}
%	The time $ \delayn $ embeds both the communication delay,
%	due to channel constraints,
%	and the computation delay,
%%	Even though the rate $ f(n) = n $ may seem natural,
%%	other rates are possible, 
%	arising if the agents preprocess the acquired measurements.
%	In practice, $ f(n) $ is to be estimated or learned from data.
%\end{rem}

\myParagraph{\titlecap{Feedback control}}
Agent $ i $ uses the received information to compute the
state mismatches $ \meas{i}{\ell^\pm}{t} $ {relative to its} neighbors,
\begin{equation}\label{eq:meas}
	\meas{i}{\ell^\pm}{t} = 
	\begin{cases}
		\xbar{i}{t} - \xbar{i\pm\ell}{t}, & 0<i\pm\ell\le N\\
		\xbar{i}{t} - \xbar{i\pm\ell\mp N}{t}, & \mbox{otherwise},
	\end{cases}
\end{equation}
%Such mismatches are exploited to compute the 
{and} the proportional control input is {given by}
\begin{equation}\label{eq:prop-control}
	\u{P,i}{t} = -\sum_{\ell=1}^{n}k_\ell\left(\meas{i}{\ell^+}{t-\taun}+\meas{i}{\ell^-}{t-\taun}\right),
\end{equation}
where measurements are delayed according to~\cref{ass:hypothesis}.

For networks with double integrator agents,
the control input $u_i(t)$ may also include a derivative term,
%Depending on the agent dynamics, the control input $ \u{i}{t} $
%may be purely proportional or include a derivative term, such as
\begin{equation}\label{eq:control-input-PD}
	\u{i}{t} = \gvel\u{P,i}{t} - \gvel\dfrac{d\xbar{i}{t}}{dt} = \gvel\u{P,i}{t} -\gvel\dfrac{d\x{i}{t}}{dt}.
\end{equation}
The derivative term in~\eqref{eq:control-input-PD} is delay free
because it only requires measurements coming from the agent itself,
which we assume {to be} available instantaneously. 
%The latter will be defined in due time.
The proportional input can be compactly written as $ \u{P}{t} = -K\xbar{}{t-\taun}=-K\x{}{t-\taun} $.
\revision{With ring topology, the feedback gain matrix is}
\begin{equation}\label{eq:feedback-matrix}
%	\begin{array}{c}
%		K \doteq K_f + K_f^\top \\
		K = \mathrm{circ}
		\left(\sum_{\ell=1}^nk_\ell, -k_1, \dots, -k_n, 0,  \dots, 0, -k_n, \dots, -k_1\right),
%	\end{array}
\end{equation}
where $ \mathrm{circ}\left(a_1,\dots,a_n\right) $ denotes the circulant matrix in $ \Real{n\times n} $
with elements $ a_1,\dots,a_n $ in the first row.

\revisiontwo{For agents with additive stochastic disturbances
	(see Sections \ref{sec:cont-time} and~\ref{sec:disc-time}),
	we consider the following problem for each $ n $.}

\begin{prob}\label{prob:variance-minimization}
	Design the feedback gains in order to minimize the steady-state variance of the consensus error,
%	\marginpar{\tiny Added both problems to highlight the optimization variables in the two cases.}
	\blue{\begin{subequations}\label{eq:problem}
		\begin{equation}\label{eq:variance-minimization-P}
		\mbox{P control:} \qquad \argmin_{K} \; \var(K),
		\end{equation}
		\begin{equation}\label{eq:variance-minimization-PD}
		\mbox{PD control:} \qquad \argmin_{\gvel,K} \; \var(\gvel,K),
		\end{equation}
	\end{subequations}}
	where
	\begin{equation}\label{}
	\var \doteq \lim_{t\rightarrow+\infty} \mathbb{E}\left[\lVert\x{}{t}\rVert^2\right]
	\end{equation}
	and w.l.o.g. we assume $ \mathbb{E}\left[\x{}{\cdot}\right] \equiv \mathbb{E}\left[\x{}{0}\right] = 0 $.
	%	and $ \varx{x} \stackrel{!}{=} $ if the system is mean-square unstable.
\end{prob}
	%!TEX ROOT = centralized_vs_distributed.tex

\section{\titlecap{Continuous-time agent dynamics}}\label{sec:cont-time}

\begin{comment}
	\done{
	\tcb{We first examine continuous-time models with single- (\autoref{sec:cont-time-single-int-model}) and double-integrator (\autoref{sec:cont-time-double-int-model}) agent dynamics. In~\autoref{sec:cont-time-single-int-control-design} we address the minimum-variance control design problem, in~\autoref{sec:numerical-results} we conduct computational experiments, and in~\autoref{sec:cont-time-single-int-trade-off} we explicitly quantify the impact of delay on the performance of distributed and centralized control strategies.}}
\end{comment}

%\mjmargin{it would be good to be specific about the meaning of "next results". Results in this section? The rest of the paper?}
\revision{We now examine continuous-time {networks} with single- (\autoref{sec:cont-time-single-int-model}) 
and double-integrator (\autoref{sec:cont-time-double-int-model}) agent dynamics,
{derive} conditions for mean-square stability,
and {compute} the steady-state variance of {a stochastically forced system}.
These {developments are} instrumental {for the formulation of the control design problem which is} used to compare different {control} architectures. 
In the optimal control problem,
the steady-state variance determines the objective function and stability conditions represent constraints.
While we first formulate and solve the problem for continuous-time dynamics, 
our results also hold for discrete-time systems;
see~\autoref{sec:disc-time}.
Also,
all \linebreak results in this section hold for generic undirected topologies.}
			%!TEX ROOT = ../../centralized_vs_distributed.tex

\subsection{\titlecap{single integrator model}}\label{sec:cont-time-single-int-model}
\done{
\tcb{The dynamics of the $ i $th agent are described by the first-order differential equation} driven by standard Brownian noise $ \noisebar{i}{\cdot} $,
\begin{equation}\label{eq:cont-time-single-int-model}
	d\xbar{i}{t} = \u{P,i}{t} dt + d\noisebar{i}{t}. %, \qquad i = 1,\dots,N
\end{equation}
\tcb{The network {error} dynamics are}
\begin{equation}\label{eq:cont-time-single-int-formation-model}
	d\x{}{t} = - K\x{}{t-\taun}dt + d\noise{}{t},
\end{equation}
where the process noise is given by $ d\noise{}{t}\sim\gauss\left(0,\Omega\Omega^\top dt\right) $.
Exploiting \tcb{symmetry of the matrix} $ K $, \tcb{we employ the change of variables $ \x{}{t} = T\xtilde{}{t} $, with $ K = T\Lambda T^\top $, to obtain $ N $ decoupled scalar subsystems with state $ \xtilde{j}{t} $, $ j=1,\dots,N $,}
\begin{equation}\label{eq:cont-time-single-int-subsystem}
	d\xtilde{j}{t} = -\gpos_j\xtilde{j}{t-\taun}dt + d\noisetilde{j}{t},
\end{equation}
where $ \gpos_j $ is the $ j $th eigenvalue of $ K $.
%\marginpar{\tiny \red{here we may want to be more precise with the equation of $ \xtilde{1}{t} $
%\eg $ d\xtilde{1}{t} \overset{a.s.}{=} 0 $}}
\blue{The subsystem with \linebreak $ \gpos_1 = 0 $ has trivial dynamics,
\ie $ d\xtilde{1}{t} \equiv 0 $,
with initial condition $ \xtilde{1}{0} = 0 $ by construction.}
For $ j \neq 1 $,  %have positive eigenvalues $ \gpos_i $ and
subsystem~\eqref{eq:cont-time-single-int-subsystem} is a single integrator %with negative feedback
driven by standard Brownian noise.

\iffalse
The subsystem with $ \gpos_1 = 0 $ has trivial dynamics
because %its state 
the mean of $ \x{}{t} $
is not controllable \tcb{with} the proportional control law~\eqref{eq:prop-control}.
\fi
}
					%!TEX ROOT = ../../centralized_vs_distributed.tex

\done{
\myParagraph{\titlecap{stability analysis}}\label{sec:setup-cont-time-single-int-stability}
\tcb{Mean-square stability} of scalar \tcb{stochastic differential equations} of the form~\eqref{eq:cont-time-single-int-subsystem} has been \tcb{addressed in the literature.}
{We build on the classical result in~\cite{KuchlerLangevinEqs} to characterize consensus stability for the multi-agent formation.} %~\eqref{eq:cont-time-single-int-formation-model}.
%leading to the following result.
\begin{prop}[Stability of CT single integrators]\label{thm:retarded-eq-steady-state}
%	\marginpar{\tiny Multi-agent is explicit now and the conditions hold for the formation.}
	{The network error $ \x{}{t} $
	is mean-square stable} if and only if 
	\begin{equation}\label{eq:cont-time-single-int-variance-condition}
		\gpos_j \in \left(0,\dfrac{\pi}{2\taun}\right), \quad j = 2,\dots,N.
	\end{equation}
	In this case, $ \x{}{t} $ is a  Gaussian process \blue{and its steady-state variance is determined by}
	\begin{equation}\label{eq:cont-time-single-int-steady-state-variance}
		\var(K) = \sum_{j=2}^{N}\varx{\gpos_j}{I}, \quad \varx{\gpos_j}{I} = \dfrac{1+\sin(\gpos_j\taun)}{2\gpos_j\cos(\gpos_j\taun)},
	\end{equation}
	where %we make explicit the dependence of $ \var $ on $ K $
	$ \varx{\gpos_j}{I} $ is the variance of the trivial solution of~\eqref{eq:cont-time-single-int-subsystem}.
\end{prop}

\begin{proof}[Sketch of Proof]
	In view of the decoupling,
	stability of~\eqref{eq:cont-time-single-int-formation-model}
	amounts to stability of all subsystems~\eqref{eq:cont-time-single-int-subsystem}, $ j=1,\dots,N $,
	with the variances of $ \x{}{t} $ and $ \xtilde{}{t} $ being equal.
	Condition~\eqref{eq:cont-time-single-int-variance-condition} and
	expression~\eqref{eq:cont-time-single-int-steady-state-variance}
	were derived in~\cite{KuchlerLangevinEqs}.
\end{proof}

\iffalse
\begin{thm}[\!\!\cite{KuchlerLangevinEqs}]
	The mean vector trivial solution of~\eqref{eq:cont-time-single-int-subsystem},
	\tcb{where $ \noisetilde{}{\cdot} $ is} standard Brownian noise,
	is mean-square stable if and only if 
	\begin{equation}\label{eq:cont-time-single-int-variance-condition}
	\gpos_j \in \left(0,\dfrac{\pi}{2\taun}\right)
	\end{equation}
	In this case, $ \xtilde{j}{t} $ is a zero-mean Gaussian process \tcb{and its steady-state variance is determined by}
	\begin{equation}\label{eq:cont-time-single-int-steady-state-variance}
	\varx{\gpos_j} = \dfrac{1+\sin(\gpos_j\taun)}{2\gpos_j\cos(\gpos_j\taun)}.
	\end{equation}
\end{thm}
\fi

While the variance of delay-free systems is bounded for any %positive semi-definite feedback matrix gain $ K $
positive eigenvalues $ \gpos_2,\dots,\gpos_N $,
%and vanishes when these go to infinity,
the presence of delay constrains a stabilizing control \tcb{according to}~\eqref{eq:cont-time-single-int-variance-condition}.
In fact,
longer delays $ \taun $ induce smaller upper bounds on the eigenvalues.

{The following result will turn useful in the control design.}
%\marginpar{\tiny This corollary is not meaningful now but is used later to assess convexity of the optimization problem and ifor the quadratic approximation.
%Also, it is referred to in Fig. 2 and in the model approximation of double integrator.}
\begin{cor}\label{lem:optimal-variance-explicit}
	\tcb{Let $ \gpos $ satisfy~\eqref{eq:cont-time-single-int-variance-condition}.
		Then} the function $ \varx{\gpos}{I} $ is strictly convex and {the minimizer} $ \opteig $ is {determined} by
	\begin{equation}\label{eq:optimal-variance-closed-form}
	\opteig = \frac{\beta^*}{\taun}, \qquad \beta^* = \cos\beta^*. %, \quad \varx{\opteig}{I} \doteq \dfrac{1+\sin(\beta^*)}{2\beta^*\cos(\beta^*)}\taun
	\end{equation}
%	where $ \beta^* \in \left(0,\nicefrac{\pi}{2}\right) $ is the unique solution of $ \beta = \cos\beta $.
\end{cor}
\begin{proof}
	Follows from standard computations over the derivatives of $ \varx{\cdot}{I} $.
	\if0\mode
	See technical report~\cite{2021arXiv210900359B}.
	\else
	See Appendices A-B in the technical report~\cite{2021arXiv210110394B}.
	\fi
\end{proof}
}
			%!TEX ROOT = centralized_vs_distributed.tex

\subsection{\titlecap{double integrator model}}\label{sec:cont-time-double-int-model}

\done{
We now \tcb{examine networks in which each agent obeys} a second-order dynamics
with the PD control input~\eqref{eq:control-input-PD}:
\begin{equation}\label{eq:cont-time-double-int-model}
	\dfrac{d^2\xbar{i}{t}}{dt^2} = \u{i}{t} + \dfrac{d\noisebar{i}{t}}{dt}.
\end{equation}
%where the $ \u{i}{t} $ is given by.
%\begin{equation}\label{eq:prop-der-control}
%	\u{i}{t} = %-\gvel\dfrac{d\xbar{i}{t}}{dt} - \gvel\u{p,i}{t} = 
%	-\gvel\dfrac{d\x{i}{t}}{dt} + \gvel\u{p,i}{t}
%	\bar{\gvel}\sum_{\ell=1}^{n}\bar{k}_\ell\left[\left(\meas{i}{\ell^+}{t-\delayn}+\meas{i}{\ell^-}{t-\delayn}\right)\right]
%\end{equation}
%where a dependence on theis added. % and  \gvel $ is to be designed.
\tcb{For simplicity,} we normalize the delay by rescaling~\eqref{eq:cont-time-double-int-model}, %follows:
\begin{equation}\label{eq:substitutions-4-normalization}
	\xbar{i}{\cdot} \leftarrow \xbar{i}{\taun\,\cdot}, \ \gvel \leftarrow \taun\gvel, \
	k_\ell \leftarrow \taun k_\ell, \ \noisebar{i}{\cdot} \leftarrow \taun \noisebar{i}{\cdot},
\end{equation}
%\cref{eq:substitutions-4-normalization} shows that any delay system
%can be brought to unit-delay
%by suitably scaling the coefficients.
%In the following, we assume the above normalization be performed in the first place.
%By re-defining the aggregate formation error $ \x{}{t} $ to be
%\begin{equation}\label{eq:multi-agent-state}
%\x{}{t} = \left[\x{1}{t}, \dots, \x{N}{t}, \dfrac{d\x{1}{t}}{dt}, \dots, \dfrac{d\x{N}{t}}{dt}\right]^\top
%\end{equation}
Stacking %the aggregate formation vector
the agent errors and their derivatives
in the formation vector, %$ \nicefrac{d\x{i}{t}}{dt} $
the error dynamics can be
%written in compact form as
%\begin{gather}
%	d\x{}{t} =
%		\left(A_0\x{}{t} + A_1\x{}{t-1}\right)dt + 
%		\begin{bmatrix}
%			0\\
%			I
%		\end{bmatrix}d\noise{}{t}\label{eq:multi-agent-state-space} \\
%	\nonumber
%	A_0 \doteq \begin{bmatrix}
%		0 & I\\
%		0 & -\gvel I
%	\end{bmatrix}, \
%	A_1 \doteq \begin{bmatrix}
%		0 & 0\\
%		-\gvel K & 0
%	\end{bmatrix}
%\end{gather}
%and $ \noise{}{\cdot} $ is composed of the $ N $ independent coordinates $ \noise{i}{\cdot} $.
%System~\eqref{eq:multi-agent-state-space} can be
decoupled as before,
%through the change of basis $ \x{}{t} = (T\otimes I_2)\xtilde{}{t} $,
%yields the following dynamics:
%\begin{align}\label{eq:multi-agent-state-space-1}
%	\begin{split}
%	d\xtilde{}{t} =
%	\left(A_0\xtilde{}{t} + \tilde{A}_1\xtilde{}{t-1}\right)dt + 
%	\begin{bmatrix}
%		0\\
%		T^\top
%	\end{bmatrix}d\noise{}{t}\\
%	\tilde{A}_1 \doteq \begin{bmatrix}
%		0 & 0\\
%		-\gvel \Lambda & 0
%	\end{bmatrix}
%	\end{split}
%\end{align}
%which can be decoupled as follows:
%where the $ i $-th subsystem is
yielding %subsystems of the form
\begin{equation}\label{eq:agent-dynamics-1}
	\dfrac{d^2\xtilde{j}{t}}{dt^2} =  -\gvel\dfrac{d\xtilde{j}{t}}{dt} - \gvel\gpos_j\xtilde{j}{t-1} + \dfrac{d\noisetilde{j}{t}}{dt}.
\end{equation}
%where $ \gpos $ is an eigenvalue of $ K $.
}
					%!TEX ROOT = ../../centralized_vs_distributed.tex

\done{
\myParagraph{\titlecap{stability analysis}}\label{sec:stability-analysis}
We have the following result.
%where $ \x{}{t} $ is the state of the system at time $ t $,
%and $ \noise{}{t} $ is a standard Brownian motion with differential $ d\noisebar{}{t}\sim\gauss(0,dt) $.
%For the sake of simplicity, 
%we normalize the delay
%through the substitutions
%\begin{equation}\label{eq:substitutions-4-normalization}
%	\x{}{t} = \xbar{}{t\tau}, \quad \gvel = \tau\bar{\gvel}, \quad \gpos = \tau\bar{\gpos}, \quad \noise{}{\cdot} = \tau \noisebar{}{\cdot}
%\end{equation}
%which yield the unit-delay equation
%\begin{equation}\label{eq:2nd-order-diff-eq-normalized}
%	\dfrac{d^2\x{}{t}}{dt^2} + \gvel\dfrac{d\x{}{t}}{dt} + \gvel\gpos\x{}{t-1} = \dfrac{d\noise{}{t}}{dt}
%\end{equation}
%The stability of~\eqref{eq:agent-dynamics-1} is characterized as follows.
\begin{prop}[Stability of CT double integrators]\label{prop:cont-time-double-int-stability}
%	Let $ (\gvel,\gpos_i)\in\Realp{2} $, $ j=1,\dots,N $.
	{The network error $ \x{}{t} $ is mean-square stable} if %, for $ j = 2,\dots,N $,
	\begin{equation}\label{eq:cont-time-double-int-stability-condition}
		\gpos_j\in\left(0,\dfrac{\beta}{\sin\beta}\right), \ 
		\gvel = \beta\tan\beta, \ 
		\beta \in \left(0,\frac{\pi}{2}\right), \ j = 2,\dots,N.
	\end{equation}
	Condition~\eqref{eq:cont-time-double-int-stability-condition} can be equivalently written as
	\begin{equation}\label{eq:cont-time-double-int-stability-region}
		\left(\gvel,\gpos_j\right) \in \setstable \doteq
		\left\lbrace (\gvel,\gpos_j)\in\Realp{2} : \gpos_j < \phi(\gvel) \right\rbrace, \ j = 2,\dots,N,
	\end{equation}
	where the implicit function $ \phi(\cdot) $ is concave increasing and
	\begin{equation}\label{eq:cont-time-double-int-phi-description}
		\phi(0) = 1, \quad \lim_{\gvel\rightarrow+\infty}\phi(\gvel) = \dfrac{\pi}{2}.
	\end{equation}
	If $ \exists j\neq1 : (\gvel,\gpos_j) \notin \overline{\setstable} $,
%	where $ \overline{\setstable} $ is the closure of $ \setstable $,
	the system is mean-square unstable.
\end{prop}
\begin{proof}
	\revisiontwo{The proof is based on~\cite{BAPTISTINI1997259}.}
	See~\cref{app:cont-time-double-int-stability}.
\end{proof}
\begin{figure}
	\centering
	\includegraphics[width=.72\linewidth]{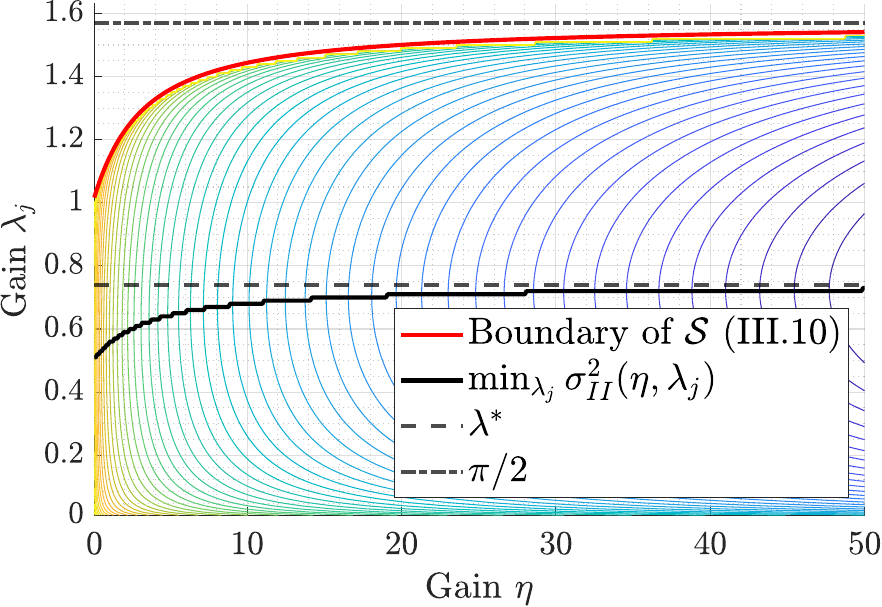}
	\caption{Level curves of the steady-state variance % $ \varx{\gvel,\gpos} $
		for the continuous-time double integrator~\eqref{eq:agent-dynamics-1}
		and points of minimum with fixed derivative gain.
	}
	\label{fig:cont-time-double-int-stab-region}
\end{figure}
\revisiontwo{
	\begin{rem}[Non-normalized delay]
		Under the original delay $ \taun $ in~\eqref{eq:cont-time-double-int-model},
		for $ j = 2,\dots,N $
		condition~\eqref{eq:cont-time-double-int-stability-condition} becomes
		\begin{equation}\label{eq:cont-time-double-int-stability-condition-rewritten}
			\gpos_j\in\left(0,\dfrac{\beta}{\taun\sin\beta}\right), \ 
			\gvel = \dfrac{\beta\tan\beta}{\taun}, \ 
			\beta \in \left(0,\frac{\pi}{2}\right).
		\end{equation}
	\end{rem}
}
%\begin{rem}[Stability conditions for feedback gains]
%While positive feedback gains ensure mean-square
%asymptotic stability for delay-free second-order systems,
%\cref{prop:cont-time-double-int-stability} states that a more restrictive condition
%applies in the presence of delay,
%similarly to the single-integrator case.
\tcb{Similar} to the single-integrator case,
\cref{prop:cont-time-double-int-stability} states that the presence of delay
requires more restrictive conditions
than positive gains.
In words, the system %~\eqref{eq:cont-time-double-int-model}
is stable
if the instantaneous \tcb{component of the}
control input in~\eqref{eq:control-input-PD}
%	(driven by the measured velocity)
is sufficiently \tcb{``strong''} compared to the delayed one.
%	which brings instability if the associated gain is too large
%\end{rem}
%When~\eqref{eq:agent-dynamics-1} is asymptotically stable, % for all $ j $,
The steady-state variance of $ \xtilde{j}{t} $ for $ j\neq1 $ can be computed
\tcb{using}~\cite[Section 4]{wangBoundedness},
\begin{equation}\label{eq:2nd-order-cont-ss-variance}
%	\lim_{t\rightarrow+\infty} \mathbb{E}[x^2(t)]
	\varx{\gvel,\gpos_j}{II} = \dfrac{1}{2\pi}\int_{-\infty}^{+\infty} \dfrac{d\omega}{|-\omega^2 + j\gvel\omega + \gvel\gpos_j\e^{-jw}|^2},
\end{equation}
and $ \var = \var(\gvel,K) = \sum_{j=2}^N \varx{\gvel,\gpos_j}{II}  $.
%which depends on both $ K $ and $ \gvel $.
\tcb{A graphical illustration of the level curves of $ \varx{\gvel,\gpos_j}{II} $ is provided in~\autoref{fig:cont-time-double-int-stab-region}.}

%\autoref{fig:cont-time-double-int-stab-region} shows
%the variance $ \varx{x} $ inside the stability region $ \mathcal{S} $.
}
					%!TEX ROOT = ../../centralized_vs_distributed.tex

%\mjmargin{blue color here is lighter than in the rest of the text where changes relative to the original submission have been made}
\done{
\myParagraph{Model Approximation}\label{sec:time-scale-separation}
\revisiontwo{Because embedding integral~\eqref{eq:2nd-order-cont-ss-variance} into an optimization problem is computationally challenging, 
	we provide an alternative tractable formulation that can be used to achieve insight into fundamental performance trade-offs.}
\tcb{As shown in~\cref{app:time-scale-separation},
	when the feedback gain $ \gvel $ is sufficiently high,
	separation of time scales\revisiontwo{~\cite{khalil2002nonlinear}} allows us to approximate~\eqref{eq:agent-dynamics-1} with first-order dynamics,}
\begin{equation}\label{eq:x-dynamics-1st-order-approximation}
	d\xtilde{j}{t} = -\gpos_j\xtilde{j}{t-1}dt + dn(t),
\end{equation}
where the \tcb{variance of Brownian motion $ n(t) $ is} inversely proportional to $ \gvel $.
In words, 
\tcb{when the damping is high enough, the derivative of $\xtilde{j}{t}$ converges to zero much faster than $\xtilde{j}{t}$, 
	which represents the dominant component of the dynamics.
Utility of this approximation is illustrated} in~\autoref{fig:cont-time-double-int-stab-region}:
with fixed $ \bar{\gvel} $, the point of minimum
%$ \argmin_{\gpos_j}\varx{\bar{\gvel},\gpos_j}{II} $
of the corresponding \tcb{1D variance} curve, \ie $ \argmin_{\gpos_j}\varx{\bar{\gvel},\gpos_j}{II} $ (solid black line),
approaches the minimizer $ \opteig $ of the single integrator \tcb{model}
(dashed black, see~\cref{lem:optimal-variance-explicit}) with increase of $ \bar{\gvel} $.
\tcb{We also note that} the variance decreases with $ \gvel $.
}
	%!TEX ROOT = ../../centralized_vs_distributed.tex

\section{\titlecap{control design}}\label{sec:cont-time-single-int-control-design}

\done{
	\myParagraph{\titlecap{single integrator model}}
	\tcb{For system~\eqref{eq:cont-time-single-int-formation-model}
		\cref{prob:variance-minimization} amounts to}
	\begin{equation}\label{eq:cont-time-single-int-variance-minimization}
		k_1^*,\dots,k_n^* = \argmin_{\{k_\ell\}_{\ell=1}^n} \; \var(K),
	\end{equation}
	\tcb{and parameterization~\eqref{eq:cont-time-single-int-subsystem} allows to rewrite it as}
	\begin{equation}\label{eq:cont-time-single-int-variance-minimization-decoupled}
		k_1^*,\dots,k_n^* = \argmin_{\{k_\ell\}_{\ell=1}^n} \; \sum_{j = 2}^N \varx{\gpos_j}{I},
	\end{equation}
	%where the circulant matrix $ K $ is uniquely defined by its eigenvalues~\cite{circulant}
	with stability condition \tcb{given by}~\eqref{eq:cont-time-single-int-variance-condition}.
	%where $ \varx{\gpos} $ is defined in~\eqref{eq:cont-time-single-int-steady-state-variance}
	%and $ \sigma(K) $ \tcb{denotes} the spectrum of $ K $.
	%$ \gpos_M \doteq \max\gpos_j < \nicefrac{\pi}{(2\taun)} $.
	\tcb{Linear dependence of the eigenvalues of $ K $ on the feedback gains~\cite{circulant}
		and~\cref{lem:optimal-variance-explicit} guarantee convexity of optimization problem~\eqref{eq:cont-time-single-int-variance-minimization-decoupled}.}
	Thus, the optimal \tcb{feedback} gains can be computed efficiently.
	
	\begin{figure}
		\centering
		\includegraphics[width=.6\linewidth]{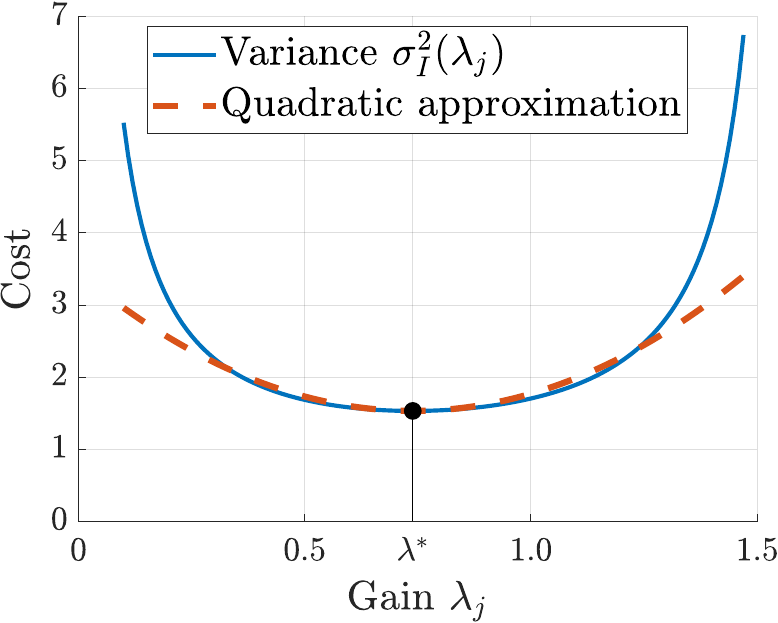}
		\caption{Exact variance function~\eqref{eq:cont-time-single-int-steady-state-variance}
			%		used in the optimal design~\eqref{eq:cont-time-single-int-variance-minimization-decoupled}
			and its quadratic approximation. % used in~\eqref{eq:quadratic-approximation}.
		}
		\label{fig:cost-comparison}
	\end{figure}
	
	\tcb{To make analytical progress and gain intuition, we also consider the following approximation of~\eqref{eq:cont-time-single-int-variance-minimization-decoupled},}
	%\marginpar{\vspace{2cm}\tiny I tried to push this figure upwards so that it's closer to the approximated optimization problem. Hope this is enough to make the whole thing more readable.}
	\begin{equation}\label{eq:quadratic-approximation}
		\tilde{k}_1^*,\dots,\tilde{k}_n^* = \argmin_{\{k_\ell\}_{\ell=1}^n} \; \sum_{j=2}^N \left(\gpos_j-\opteig\right)^2,
	\end{equation}
	\tcb{which} squeezes the spectrum of $ K $ about the ``optimal" \tcb{eigenvalue $\opteig$. The variance $ \varx{\cdot}{I} $ can be approximated with a quadratic function around its minimum because it is strictly convex, differentiable in the stability region, and it blows up at the boundaries $ \{0,\nicefrac{\pi}{2}\} $, see~\autoref{fig:cost-comparison}.}
	\begin{prop}[Near-optimal proportional control]\label{prop:subopt-gain}
		The solution of \tcb{problem~\eqref{eq:quadratic-approximation} is determined by 
			\[
			\tilde{k}_\ell^* 
			\; \equiv \; 
			\tilde{k}^*
			\; \doteq \;
			\dfrac{\opteig}{2n+1}.
			\]}
	\end{prop}
	\begin{proof}
		The result follows by applying properties of the DFT to~\eqref{eq:quadratic-approximation}.
		\if0\mode
		See technical report~\cite{2021arXiv210900359B}.
		\else
		See Appendices C-D in the technical report~\cite{2021arXiv210110394B}.
		\fi
	\end{proof}
	\cref{prop:subopt-gain} \tcb{shows that spatially-constant feedback gains provide good performance even when spatially-varying feedback gains are allowed. According to~\cref{lem:optimal-variance-explicit}, the suboptimal gain $ \tilde{k}^*$ decreases with the delay $ \taun $ and with the number of agents involved in the feedback loops, thereby reflecting benefits of communication.
	}
}
			%!TEX ROOT = ../centralized_vs_distributed.tex

\begin{figure*}%
	\centering
	\subfloat[Continuous-time single integrator.% with $ \taun = 0.1n $.
	]{
		\centering
		\includegraphics[width=.22\linewidth]{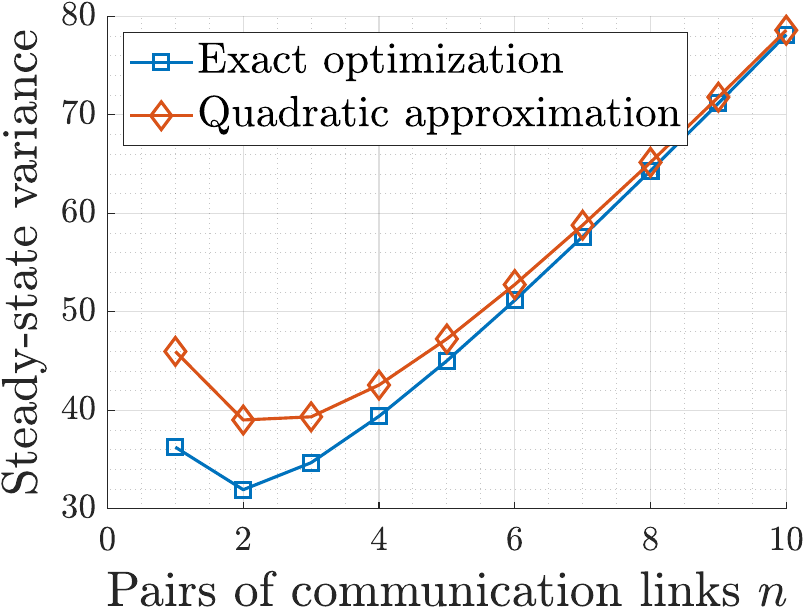}
		\label{fig:cont-time-single-int-opt-var}
	}%
	\hfil
	\subfloat[Continuous-time double integrator.% with $ \taun = 0.1n $.
	]{
		\centering
		\includegraphics[width=.22\linewidth]{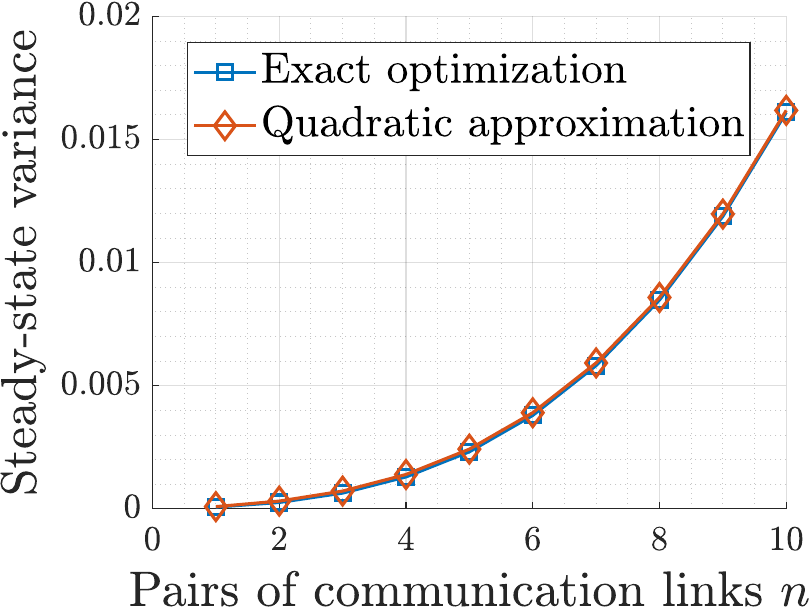}
		\label{fig:cont-time-double-int-opt-var}
	}%
	\hfil
	\subfloat[Discrete-time single integrator.% with $ \taun = n $.
	]{
		\centering
		\includegraphics[width=.22\linewidth]{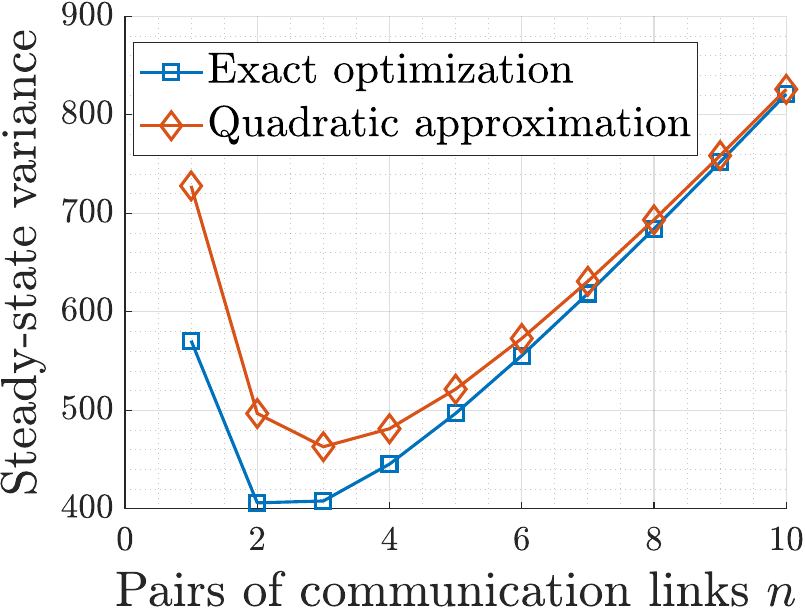}
		\label{fig:disc-time-single-int-opt-var}
	}%
	\hfil
	\subfloat[Discrete-time double integrator.% with $ \taun = n $.
	]{
		\centering
		\includegraphics[width=.22\linewidth]{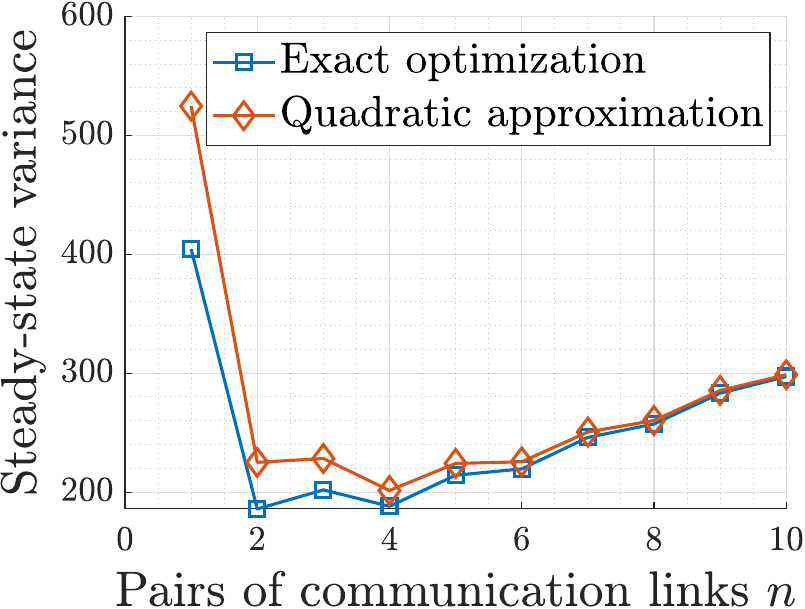}
		\label{fig:disc-time-double-int-opt-var}
	}%
	\caption{Optimal
		and suboptimal steady-state scalar variances with linear delay increase
		for different agent dynamics.
		%		The used delay increase rates are $ \taun = 0.1n $ for continuous-time and $ \taun = n $ for discrete-time models.
		%		All curves exhibit nontrivial trade-offs over the number of links. %, with the optimal amounts of links smaller than the maximum available.
	}
	\label{fig:opt-var}
\end{figure*}
					%!TEX ROOT = centralized_vs_distributed.tex

\done{
\myParagraph{\titlecap{double integrator model}}
%As for the double-integrator agent dynamics,
\tcb{Approximation~\eqref{eq:x-dynamics-1st-order-approximation}} and~\autoref{fig:cont-time-double-int-stab-region}
show that,
for sufficiently large $ \gvel $,
the variance of the double-integrator subsystem~\eqref{eq:agent-dynamics-1}
	%vanishes for \tcb{infinitely large derivative feedback gain, $ \gvel=+\infty $.
%	The same approximation implies that,	
%	the variance of the 
	has structure similar to the single integrator, \ie $ \varx{\gvel,\gpos_j}{II} \approx c \varx{\gpos_j}{I} $ for some ``small" $ c > 0 $.
%	the proportional gains should be designed as for the single-integrator case.
%\mjmargin{I don't understand the red equation. Perhaps a better way to write it is given below?}
Thus, we approximate the control design~\eqref{eq:variance-minimization-PD} as
\begin{equation}\label{eq:cont-time-double-int-min-var-simplified}
%	\argmin_{\gvel,\{k_\ell\}_{\ell=1}^n} \; \lim_{t\rightarrow+\infty} \mathbb{E}\left[\lVert\x{}{t}\rVert^2\right] \ \approx \
%\argmin_{\gvel,\{k_\ell\}_{\ell=1}^n} \; \sum_{j=2}^N \varx{\gvel,\gpos}{II} \ \approx \
\tilde{\gvel}^*, \, \argmin_{\{k_\ell\}_{\ell=1}^n} \; \sum_{j=2}^N \varx{\gpos_j}{I},
\end{equation}
where $ \tilde{\gvel}^* $ is chosen beforehand so that
the time-scale separation \tcb{argument provides} a reasonable approximation~\eqref{eq:x-dynamics-1st-order-approximation}.
%approximation~\eqref{eq:x-dynamics-1st-order} is good enough.
In particular, the \tcb{optimization problem for proportional feedback gains in~\eqref{eq:cont-time-double-int-min-var-simplified}}
coincides with the control design for single integrators~\eqref{eq:cont-time-single-int-variance-minimization-decoupled},
with the exception that the stability condition is now \tcb{given by $ \gpos_j < \phi(\tilde{\gvel}^*) $,
	$ j = 2,\dots,N $; see~\eqref{eq:cont-time-double-int-stability-region}.}
%In particular, the feasibility constraint tends to $ \gpos_M < \nicefrac{\pi}{2} $ as $ \tilde{\gvel}^* $ grows.
%namely the constraint for single integrators.

\iffalse
\canOmit{Alternatively, letting $ \tilde{\gvel}^* = +\infty $ at first,
one can choose the single-integrator optimal gains $ \{k_\ell\}_{\ell=1}^n $
according to~\eqref{eq:cont-time-single-int-variance-minimization-decoupled},
and then select a finite value $ \tilde{\gvel}^* > \phi^{-1}(\gpos_M^*) $ to guarantee stability,
where $ \gpos_M^* $ is the optimal spectral radius of $ K $.}
\setlength\marginparwidth{25pt}
\marginpar{\vspace{-4em} \small can \\ remove}
\fi

\begin{rem}[Convexity enables comparison]
	\revision{Convexity of the optimal control design
		problems~\eqref{eq:cont-time-single-int-variance-minimization-decoupled}--\eqref{eq:cont-time-double-int-min-var-simplified}
		enables both efficient numerical computations 
		of the optimal feedback gains for \textit{given} $ n $
		and fair comparison of the best achievable performance for \textit{different} values of $ n $.}
\end{rem}

\begin{rem}[Gain scaling]
	The \tcb{optimal feedback} gains $ \{k_\ell^*\}_{\ell=1}^n $ and $ \tilde{\eta}^* $ 
	are to be scaled \tcb{by} $ \nicefrac{1}{\taun}$
	according to~\eqref{eq:substitutions-4-normalization}.
	%In particular, long delays induce small gains.
\end{rem}

\begin{rem}[Optimal design for double integrators]
	%	Given any value of $ \gvel $,
	%	Problem~\eqref{eq:cont-time-double-int-min-var-simplified}
	\tcb{Local minimizer of the original problem approximated by~\eqref{eq:cont-time-double-int-min-var-simplified}
		can be solved using the gradient-based method proposed in~\cite{8358743}.
		However, 
		this approach has no guarantees of global optimality
		and its computational %further, even if it has polynomial
		complexity is impractical for large-scale systems.
		In contrast,
		convex approximation~\eqref{eq:cont-time-double-int-min-var-simplified} draws a parallel 
		to the optimal design for the single-integrator model and provides insight into a \tradeoff.} % embedded in such networked systems.
\end{rem}
}
			%!TEX ROOT = centralized_vs_distributed.tex

\done{
\subsection{\titlecap{General symmetric network topology}}\label{sec:generic-topology}
Even though we \revision{utilized ring topology to derive analytical results (see~\autoref{sec:cont-time-single-int-trade-off})},
the control design can be extended to general undirected networks with symmetric feedback gain matrices $ K $.
For the single integrator model, this reads
\begin{equation}\label{eq:generic-problem}
K^* = \argmin_K \; \var(K).
\end{equation}
\tcb{The steady-state} {network error} variance $ \var(K) $ is a convex \tcb{function}
if and only if $ \varx{\gpos_j}{I} $ is convex~\cite{davis1957all},
which is proved in~\cref{lem:optimal-variance-explicit} for continuous-time
and \tcb{in~\cref{app:disc-time-single-int-variance-explicit} for} discrete-time systems. %and the system eigenvalues
%and is thus convex in the feedback gains
The optimal gains can then be found numerically via gradient-based methods, where \tcb{gradients of the eigenvalues can be computed using analytical~\cite{doi:10.2514/3.7211,doi:10.2514/2.1119} or numerical~\cite{10.1115/1.2888195} methods.}
On the other hand,
the derivative {feedback gain in $ \varx{\gvel,\gpos_j}{II} $ prevents us from establishing convexity for second-order systems in general.}
However, if $ \varx{\gvel,\gpos_j}{II} $ is convex in each coordinate%
\footnote{This can be checked for discrete-time double integrators, see~\cref{app:disc-time-single-int-variance-explicit}.},
the design \tcb{problem} can be solved by alternatively optimizing proportional and derivative gains \tcb{and the \tradeoff can be studied irrespective of the particular topology.}}
	%!TEX ROOT = ../../centralized_vs_distributed.tex

\section{{\titlecap{the centralized-distributed trade-off}}}\label{sec:numerical-results}

\revision{In the previous sections we formulated the optimal control problem for a given controller architecture
(\ie the number of links) parametrized by $ n $
and showed how to compute minimum-variance objective function and the corresponding constraints.
In this section, we present our main result:
%\red{for a ring topology with multiple options for the parameter $ n $},
we solve the optimal control problem for each $ n $ and compare the best achievable closed-loop performance with different control architectures.\footnote{
\revision{Recall that small (large) values of $ n $ mean sparse (dense) architectures.}}
For delays that increase linearly with $n$,
\ie $ f(n) \propto n $, 
we demonstrate that distributed controllers with} {few communication links outperform controllers with larger number of communication links.}

\textcolor{subsectioncolor}{Figure~\ref{fig:cont-time-single-int-opt-var}} shows the steady-state variances
obtained with single-integrator dynamics~\eqref{eq:cont-time-single-int-variance-minimization}
%where we compare the standard multi-parameter design 
%with a simplified version \tcb{that utilizes spatially-constant feedback gains
and the quadratic approximation~\eqref{eq:quadratic-approximation} for \revision{ring topology}
with $ N = 50 $ nodes. % and $ n\in\{1,\dots,10\} $.
%with $ N = 50 $, $ f(n) = n $ and $ \tau_{\textit{min}} = 0.1 $.
%\autoref{fig:cont-time-single-int-err} shows the relative error, defined as
%\begin{equation}\label{eq:relative-error}
%	e \doteq \dfrac{\optvarx-\optvar}{\optvar}
%\end{equation}
%where $ \optvar $ and $ \optvarx $ denote the the optimal and sub-optimal scalar variances, respectively.
%The performance gap is small
%and becomes negligible for large $ n $.
{The best performance is achieved for a sparse architecture with  $ n = 2 $ 
in which each agent communicates with the two closest pairs of neighboring nodes. 
This should be compared and contrasted to nearest-neighbor and all-to-all 
communication topologies which induce higher closed-loop variances. 
Thus, 
the advantage of introducing additional communication links diminishes 
beyond}
{a certain threshold because of communication delays.}

%For a linear increase in the delay,
\textcolor{subsectioncolor}{Figure~\ref{fig:cont-time-double-int-opt-var}} shows that the use of approximation~\eqref{eq:cont-time-double-int-min-var-simplified} with $ \tilde{\gvel}^* = 70 $
identifies nearest-neighbor information exchange as the {near-optimal} architecture for a double-integrator model
with ring topology. 
This can be explained by noting that the variance of the process noise $ n(t) $
in the reduced model~\eqref{eq:x-dynamics-1st-order-approximation}
is proportional to $ \nicefrac{1}{\gvel} $ and thereby to $ \taun $,
according to~\eqref{eq:substitutions-4-normalization},
making the variance scale with the delay.

%\mjmargin{i feel that we need to comment about different results that we obtained for CT and DT double-intergrator dynamics (monotonic deterioration of performance for the former and oscillations for the latter)}
\revision{\textcolor{subsectioncolor}{Figures~\ref{fig:disc-time-single-int-opt-var}--\ref{fig:disc-time-double-int-opt-var}}
show the results obtained by solving the optimal control problem for discrete-time dynamics.
%which exhibit similar trade-offs.
The oscillations about the minimum in~\autoref{fig:disc-time-double-int-opt-var}
are compatible with the investigated \tradeoff~\eqref{eq:trade-off}:
in general, 
the sum of two monotone functions does not have a unique local minimum.
Details about discrete-time systems are deferred to~\autoref{sec:disc-time}.
Interestingly,
double integrators with continuous- (\autoref{fig:cont-time-double-int-opt-var}) ad discrete-time (\autoref{fig:disc-time-double-int-opt-var}) dynamics
exhibits very different trade-off curves,
whereby performance monotonically deteriorates for the former and oscillates for the latter.
While a clear interpretation is difficult because there is no explicit expression of the variance as a function of $ n $,
one possible explanation might be the first-order approximation used to compute gains in the continuous-time case.
%which reinforce our thesis exposed in~\autoref{sec:contribution}.

%\begin{figure}
%	\centering
%	\includegraphics[width=.6\linewidth]{cont-time-double-int-opt-var-n}
%	\caption{Steady-state scalar variance for continuous-time double integrators with $ \taun = 0.1n $.
%		Here, the \tradeoff is optimized by nearest-neighbor interaction.
%	}
%	\label{fig:cont-time-double-int-opt-var-lin}
%\end{figure}
}

\begin{figure}
	\centering
	\begin{minipage}[l]{.5\linewidth}
		\centering
		\includegraphics[width=\linewidth]{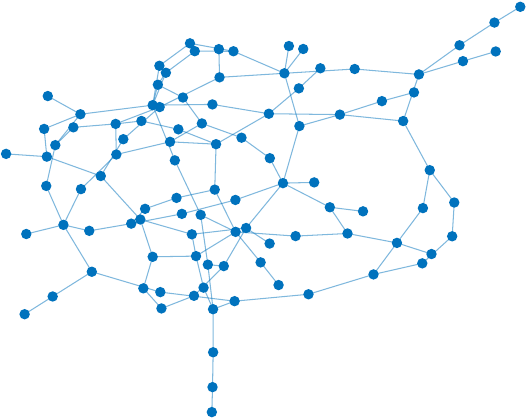}
	\end{minipage}%
	\begin{minipage}[r]{.5\linewidth}
		\centering
		\includegraphics[width=\linewidth]{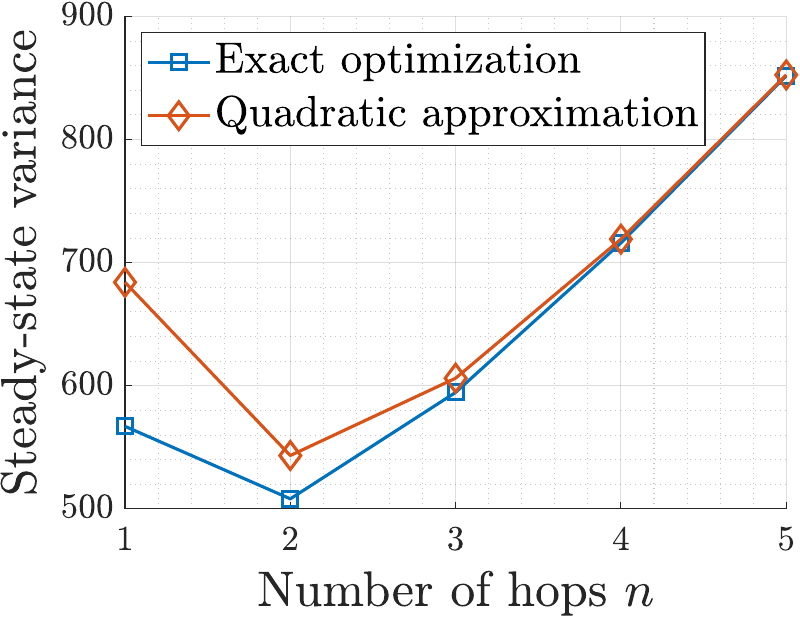}
	\end{minipage}
	\caption{Network topology and its optimal {closed-loop} variance.}
	\label{fig:general-graph}
\end{figure}

Finally,
\autoref{fig:general-graph} shows the optimization results for a random graph topology with discrete-time single integrator agents. % with a linear increase in the delay, $ \taun = n $.
Here, $ n $ denotes the number of communication hops in the ``original" network, shown in~\autoref{fig:general-graph}:
as $ n $ increases, each agent can first communicate with its nearest neighbors,
then with its neighbors' neighbors, and so on. For a control architecture that utilizes different feedback gains for each communication link
	(\ie we only require $ K = K^\top $) we demonstrate that, in this case, two communication hops provide optimal closed-loop performance. % of the system.}

Additional computational experiments performed with different rates $ f(\cdot) $ show that the optimal number of links increases for slower rates: 
for example, 
the optimal number of links is larger for $ f(n) = \sqrt{n} $ than for $ f(n) = n $. 
\revision{These results are not reported because of space limitations.}
			%!TEX ROOT = ../../centralized_vs_distributed.tex

\subsection{{\titlecap{ring topology: analytical insight into the trade-off}}}\label{sec:cont-time-single-int-trade-off}

\done{
%In particular, %as the delay increases according to~\cref{ass:hypothesis},
%the steady-state error variance is minimized by a network topology
%whose optimal communication neighborhood size $ n^* $
%is smaller than the maximum number of links,
%%(corresponding to the complete graph, \ie centralized control)
%in general.
%represented by the optimal number of link pairs $ n^* $.
%Beyond such a threshold, further enlarging the feedback loops
%accrues little benefit as opposed to the penalization of the dynamics due to increased latency,
%so that even the optimal control design yields a worse performance than with fewer links.

For \revision{\tcb{a ring topology with continuous-time single-integrator agent dynamics,
	a \tradeoff can be explicitly quantified.
	By utilizing~\cref{prop:subopt-gain} to compute the feedback gains,
	the objective function can be factorized as}}
%	 \mjmargin{roman font for subscripts}
\begin{equation}\label{eq:cont-time-single-int-trade-off-mult}
	\var = \underbrace{f(n)}_{\tilde{J}_{\textrm{latency}}(n)} \cdot \	\underbrace{\sum_{j=2}^N \tilde{C}_{j}^*(n)}_{\tilde{J}_{\textrm{network}}(n)},
\end{equation}
%where the optimal variance $ \rho^* $ only depends on the delay $ \taun $.
%$ \rho^* \doteq \rho(\opteig) $ is the minimum of the variance function associated with the scalar subsystem
where %$ \tilde{K}^* $ is the suboptimal feedback matrix, %with the gains computed as per~\cref{prop:subopt-gain}
$ \sigma_{I}^{2}(\tilde{\gpos}_j^*) =  \tilde{C}_{j}^*(n)\taun $
%with $ \tilde{\gpos}_j^* $ a suboptimal eigenvalue of $ K $,
%with the suboptimal gain $ \tilde{k}^* $ and
and $ \tilde{C}_{j}^*(n) $ only depends on $ n $
and can be computed \tcb{exactly; see Appendix~\ref{app:cont-time-single-int-suboptimal-variance-computation}.} 
%\mjmargin{red sentence is a bit mouthful (e.g., $ \tilde{\gpos}^* $ are linear in $ \opteig $) -- please try to rewrite.}
{This holds because %the gain $ \tilde{k}^* $,
	the suboptimal eigenvalues can be expressed as $ \tilde{\gpos}_j^* = \tilde{c}_j^*(n)\opteig $ %
%are linear in $  $ through coefficients $ $
%that only \tcb{depend} on \tcb{$ n $
	(cf.~\cref{prop:subopt-gain}).}
%where $ \tilde{\gpos}^* = \tilde{c}^*(n)\opteig $ % is used.
%and $ \tilde{c}^*(n) $ only depends on $ n $.
Such a decomposition can be interpreted as a
%The factors $ \tilde{c}^*(n) $ and $ \opteig $
decoupling of \tcb{the impact of network ($ \tilde{c}_j^*(n) $) and latency ($ \opteig $) effects on the control design.}
% being linear . % (c.f.~\cref{prop:subopt-gain}).
%The factors $ \tilde{c}_i^* $ are easily computed from the analytical expression of the spectrum of $ K $
%as the Discrete Fourier Transform of the first row.
%Notice that~\eqref{eq:single-int-variance-minimization-rewritten} is similar to~\eqref{eq:trade-off},
%with the only exception of the additional (increasing) factor $ f(n) $ which multiplies the network-related cost.
%The observed trade-off over $ n $ arises from the opposite behaviors of the two addends in~\eqref{eq:single-int-variance-minimization-rewritten}.
%\autoref{fig:trade-off} illustrates such trends for the curves in Figs.~\ref{fig:cont-time-single-int-opt-var}--\ref{fig:disc-time-single-int-opt-var},
%where the average gap is proportional to the summation in~\eqref{eq:single-int-variance-minimization-rewritten}.
%On the one hand,
%~\autoref{fig:coeff-Ctilde} shows ,
By inspection, it can be seen that $ \tilde{J}_{\textrm{network}}(n) $ is a decreasing function of $n$
		and that $ \tilde{J}_\textrm{latency}(n) $ is determined by $ f(n) $.
		Furthermore, when $ f(\cdot) $ is sublinear,
		the above expression can be equivalently written in form~\eqref{eq:trade-off}, %according to~\eqref{eq:trade-off}:
\begin{equation}\label{eq:cont-time-single-int-trade-off}
	\var = \underbrace{f(n)\cdot\sum_{j=2}^N\left(\tilde{C}_{j}^*(n)-C^*\right)}_{\tcb{J_{\textrm{network}}(n)}} + 
	\underbrace{(N-1)C^*f(n)}_{\tcb{J_{\textrm{latency}}(n)}},
\end{equation}
where $ \varx{\opteig}{I} = C^*\taun $ is the optimal variance \tcb{according to~\eqref{eq:cont-time-single-int-steady-state-variance} and~\cref{lem:optimal-variance-explicit}.}
Indeed,
the summation decreases with superlinear rate,
so that $ J_{\textrm{network}}(n) $ is a decreasing sequence.
The terms in $ J_{\textrm{network}}(n) $,
each associated with a decoupled subsystem~\eqref{eq:cont-time-single-int-subsystem},
illustrate benefits of communication:
as $ n $ increases, the eigenvalues of $ K $ have more degrees of freedom
and can squeeze more tightly about $ \opteig $,
reducing performance gaps between subsystems and theoretical optimum.
We note that $ J_{\textrm{network}}(n) $ vanishes for the fully connected architecture.

%\marginpar{\vspace{-2cm}\tiny I reduced this comment and tried to sharpen a little bit to reinforce our thesis is valid for general systems.
%	This is not the main point of the section though, so if you feel this is not strong enough you may remove it.}
{Even though analogous expressions could not be obtained for other dynamics,
	the curves in~\autoref{fig:opt-var} exhibit trade-offs which are consistent with the above analysis.}
}
	%!TEX ROOT = ../centralized_vs_distributed.tex

\section{\titlecap{Discrete-time agent dynamics}}\label{sec:disc-time}

\blue{We now consider discrete-time agent dynamics
to illustrate that the afore-established fundamental trade-offs hold in this case as well.}
%to address realistic communication.
%As shown in~\autoref{sec:numerical-results},
%this does not affect the fundamental \tradeoff observed for continuous-time systems.
In what follows, we denote time instants by $ \{k\}_{k\in\mathbb{N}} \doteq \{kT\}_{k\in\mathbb{N}} $,
$ T $ being the sampling time.
Similarly, we re-define the delay as the number of delay steps $ \taun \doteq \ceil{\nicefrac{\taun}{T_s}} $.
					%!TEX ROOT = ../../centralized_vs_distributed.tex

\myParagraph{\titlecap{agent models}}
The discrete-time versions of the agent dynamics considered in~\autoref{sec:cont-time} are given by
\begin{equation}\label{eq:disc-time-single-int-model}
	\xbar{i}{k+1} = \xbar{i}{k} + \u{P,i}{k} + \noisebar{i}{k},
\end{equation}
for the single-integrator model, with $ \noisebar{i}{\cdot}\sim\gauss(0,1) $, and
%with process noise $ \noise{}{\cdot} \sim \gauss(0,\Omega\Omega^\top) $.
%the multi-agent system can be written compactly in state-space form as
%\begin{equation}\label{eq:disc-time-single-int-state-space}
%\x{}{k+1} = \x{}{k} - K\x{}{k-\taun} + \noise{}{k}
%\end{equation}
\begin{equation}\label{eq:disc-time-double-int-model}
	\begin{aligned}
		\xbar{i}{k+1} &= \xbar{i}{k} + \zbar{i}{k}\\
		\zbar{i}{k+1} &= (1-\gvel)\zbar{i}{k} + \gvel\u{P,i}{k} + \noisebar{i}{k},
	\end{aligned}
\end{equation}
for the double-integrator model,
with $ \u{P,i}{k} $ defined in~\eqref{eq:prop-control}.
					%!TEX ROOT = ../../centralized_vs_distributed.tex

\myParagraph{Stability Analysis}\label{sec:disc-time-stability-analysis}
%Writing the formation dynamics in vector-matrix form,
The formation error dynamics can be decoupled analogously to the continuous-time models.
The decoupled subsystems are asymptotically stable
%poles of the spectral factors $ W_i(z) $
%such that $ \xtilde{i}{k} = W_i(z)\noisetilde{i}{k} $.
%which generates $ \xtilde{i}{t} $ by filtering white noise.
%Consider the generic system
%\begin{equation}\label{eq:disc-time-single-int-subsystem}
%		\x{}{k+1} = \x{}{k} -\gpos\x{}{k-\delayn} + \noise{}{k}
%\end{equation}
%In~\eqref{eq:disc-time-single-int-subsystem}, $ \x{}{k} $ is generated from white noise through the filter
%\begin{equation}\label{eq:disc-time-single-int-transfer-function}
%	W(z) = \dfrac{1}{z-1+\gpos z^{-\delayn}}
%\end{equation}
%In particular, the systems are asymptotically stable
if and only if all the roots of their associated characteristic polynomials
%\ie the denominator of $ W_i(z) $,
lie inside the unit circle in the complex plane.

\begin{figure}
	\centering
	\begin{tikzpicture}[scale=1]
		\begin{axis}[xmin=0,xmax=20,ymin=0,ymax=2,
			ytick = {0,1,pi/2},yticklabels = {$ 0 $,$ 1 $,$ \dfrac{\pi}{2} $},
			tick label style={font=\normalsize},
			xlabel=$ \taun $,ylabel=$ \gpos_j\taun $,label style={font=\large},
			trig format plots=rad,
			legend pos=south east,legend style={font=\normalsize},legend cell align={left},
			width=.9\linewidth,height=.5\linewidth]
			\addplot[draw=black,pattern=north west lines,pattern color=black!40,area legend] (0,0) rectangle (20,pi/2);
			\addplot[domain=1:20,samples=20,mark=none,ycomb,ultra thick] {x*2*sin(pi/(2*(2*x+1)))};
			\addplot[draw=black,ultra thick] (20,0) -- (20,1);
			\legend{Continuous-time~\eqref{eq:cont-time-single-int-variance-condition},,Discrete-time~\eqref{eq:disc-time-single-int-stability-condition}};
		\end{axis}
	\end{tikzpicture}
	\caption{Stability regions of decoupled single integrators.}
	\label{fig:stability-region-single-int}
\end{figure}
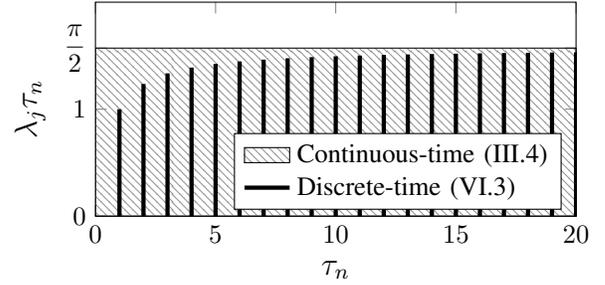

In general, given a delay $ \taun $,
stability conditions with respect to the control gains can be derived
%the location of the system poles %of~\eqref{eq:disc-time-single-int-characteristic-polinomial}--\eqref{eq:disc-time-double-int-characteristic-polinomial}
in the form of polynomial inequalities through the Jury criterion.
For the single-integrator case,
one simple condition can be computed analytically.  %by exploiting the root locus associated with $ \gpos_i $
\begin{prop}[Stability of DC single integrators]\label{prop:disc-time-single-int-stability}
%	Let $ \gpos \in\Realp{} $.
%	Then, system~\eqref{eq:disc-time-single-int-model}
	{The network error $ \x{}{t} $ is mean-square stable} if and only if
%	the eigenvalues of $ K $ satisfy
	\begin{equation}\label{eq:disc-time-single-int-stability-condition}
		\gpos_j \in \left(0,2\sin\left(\dfrac{\pi}{2}\dfrac{1}{2\taun+1}\right)\right), \quad j = 2,\dots,N.
%		\gpos < \gpos_{\textit{th}} \doteq 2\sin\left(\dfrac{\pi}{2}\dfrac{1}{2\tau+1}\right)
	\end{equation}
\end{prop}
The upper bound in~\eqref{eq:disc-time-single-int-stability-condition} approaches its continuous-time counterpart~\eqref{eq:cont-time-single-int-variance-condition} from below
as the delay steps tend to infinity (see~\autoref{fig:stability-region-single-int}). %, where $ \sin(\star) \approx \star $.
%Indeed, given the same absolute delay,
%a finer sampling yields more delay steps.
%and thus at the limit the discretized dynamics
%converges to the continuous-time one. and retrieves the same constraint.
%In general, condition~\eqref{eq:disc-time-single-int-stability-condition} is tighter than~\eqref{eq:cont-time-single-int-variance-condition}.
%On the other hand, the asymptotic behavior of the threshold gain
%suggests that the gap between continuous-time and discretized systems only matters
%when the delay is comparable with the sampling time,
%while, when the former gets too long, the loss of feedback information
%neglects the dynamics discretization.
A discussion on general stability conditions
and the proof of~\cref{prop:disc-time-single-int-stability} are provided in~\cref{app:disc-time-single-int-stability}.
The basic argument is the same as for the continuous-time case.
%\autoref{fig:stability-region-single-int} compares the stability regions
%in the $ (\taun,\gpos_j\taun) $ plane
%for continuous- and discrete-time single integrators.
%In particular, the stability region of the discrete-time systems is strictly contained in the continuous-time one.

					%!TEX ROOT = ../../centralized_vs_distributed.tex

\begin{figure}
	\centering
	\begin{minipage}{.5\linewidth}
		\centering
		\includegraphics[width=\linewidth]{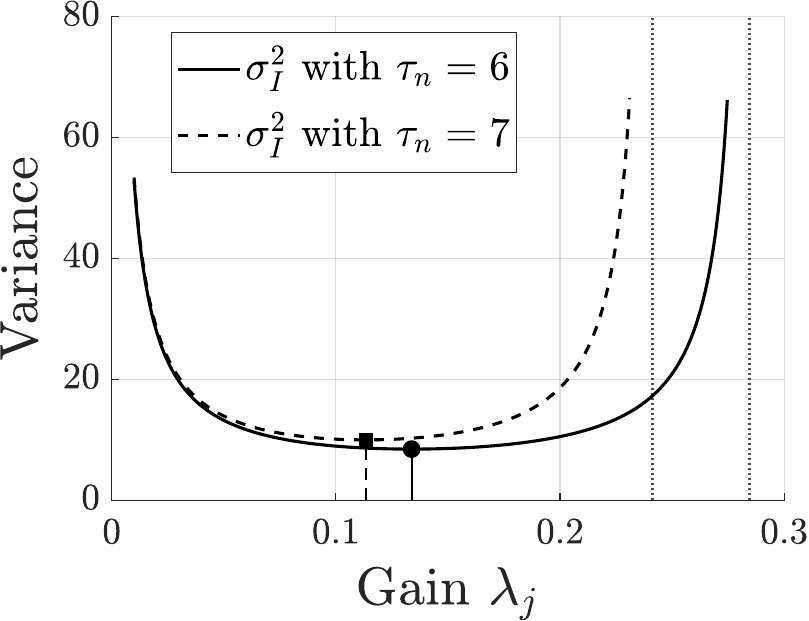}
	\end{minipage}%
	\hfil
	\begin{minipage}{.5\linewidth}
		\centering
		\includegraphics[width=\linewidth]{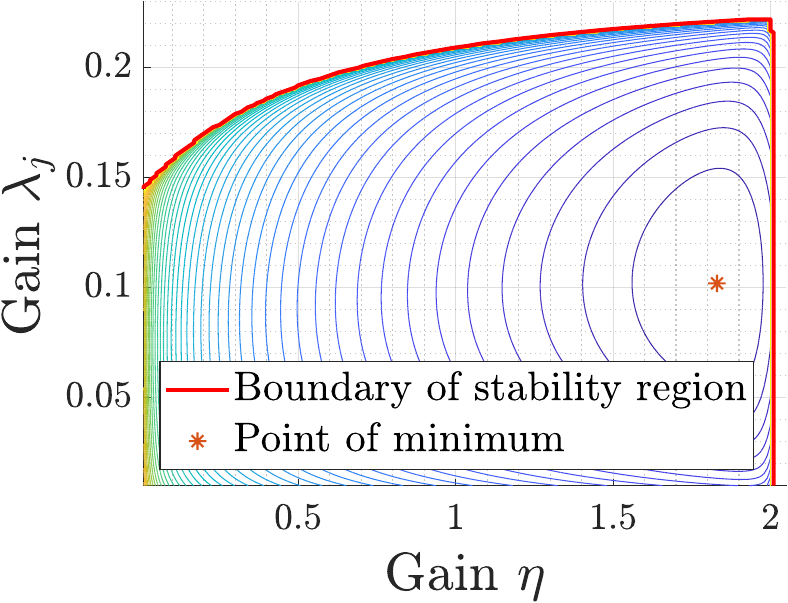}
	\end{minipage}
	\caption{Typical profiles of the steady-state variance
		for decoupled discrete-time single integrators (left) and double integrators (right). %with $ \delayn \in \{5,7\} $.
%		The minima are highlighted by markers.
	}
	\label{fig:disc-time-var}
\end{figure}

\myParagraph{\titlecap{performance evaluation}}\label{sec:disc-time-single-int-moment-matching}
With fixed parameters,
the steady-state variance of each decoupled subsystem
can be computed numerically via the Wiener–Khintchine formula. % recalled in~\autoref{app:disc-time-single-int-variance-explicit}.
Also, for any given value of $ \taun $,
%gain-parametric 
a closed-form expression of the variance
%whose convexity can be easily assessed,
can be obtained via moment matching through a recursive formula, see~\cref{app:disc-time-single-int-variance-explicit}.
Such closed-form expressions have been used for our computational experiments illustrated in~\autoref{fig:opt-var}.
\textcolor{subsectioncolor}{Figure~\ref{fig:disc-time-var}} shows the typical profiles
of the variance function for decoupled subsystems with single- and double-integrator dynamics
(see~\eqref{eq:disc-time-single-int-decoupled} and~\eqref{eq:disc-time-double-int-decoupled} in~\cref{app:disc-time-single-int-variance-explicit},
respectively).
%For the one-dimensional case (single integrator),
%convexity of the variance function $ \var{\gpos_i} $ can be checked
%by studying the second derivative. % solving a system of inequalities.
%%over the decoupled subsystems.
%\autoref{fig:disc-time-var} shows the level curves of the variance,
%which looks convex also for the double integrator.

	%!TEX ROOT = ../centralized_vs_distributed.tex

\section{\titlecap{Conclusion and future research}}\label{sec:conclusion}

{We study minimum-variance control design problem for undirected networks with both 
continuous- and discrete-time agent dynamics in the presence of communication delays.
When feedback delays increase with the number of communication links, 
we identify fundamental performance trade-offs and
show that distributed control architectures can offer superior performance to centralized ones that utilize all-to-all information exchange.
Our hope is to pave the way to a new body of research which will
enable control design with a deeper understanding of the fundamental behavior and limitations of large-scale wireless network systems.
Future work will focus on extending our results to other classes of control 
problems which include more complex system dynamics and communication models, 
more realistic information about structure of delays in a distributed scenario, 
as well as different cost functions.}
%for example, the \tradeoff may be further investigated
%under %different effects of latency in the system dynamics (\eg multiple delays),
%multiple, stochastic or time-varying delays,
%unreliable communication, or
%heterogeneous agents with different delay functions. % computational capabilities and/or transmission power.
	
	\if0\mode
	\bibliographystyle{IEEEtran}
	\bibliography{bibfile}
	\else

	\fi
	
	\appendix
	\numberwithin{equation}{subsection}
	%!TEX ROOT = ../centralized_vs_distributed.tex

\subsection{Proof of~\cref{prop:cont-time-double-int-stability}}\label{app:cont-time-double-int-stability}

%Let us consider the following generic scalar stochastic retarded linear system driven by a standard Brownian motion:
%\begin{equation}\label{eq:2nd-order-diff-eq}
%\dfrac{d^2\x{}{t}}{dt^2} + \gvel\dfrac{d\x{}{t}}{dt} + \gvel\gpos\x{}{t-1} = \dfrac{d\noise{}{t}}{dt}
%\end{equation}
The error dynamics equation with agent model~\eqref{eq:cont-time-double-int-model} reads
\begin{equation}\label{eq:multi-agent-state-space}
	\begin{array}{c}
		d\x{}{t} =
			\left(A_0\x{}{t} + A_1\x{}{t-1}\right)dt + Bd\noisebar{}{t}, \\[10pt]
		A_0 = \begin{bmatrix}
			0 & I\\
			0 & -\gvel I
		\end{bmatrix}, \
		A_1 = \begin{bmatrix}
			0 & 0\\
			-\gvel K & 0
		\end{bmatrix}, \
		B = \begin{bmatrix}
			0\\
			I
		\end{bmatrix},
	\end{array}
\end{equation}
with $ \noisebar{}{t} $ standard $ N $-dimensional Brownian motion.
The decoupling~\eqref{eq:agent-dynamics-1}
is obtained from~\eqref{eq:multi-agent-state-space} through
the change of basis $ \x{}{t} = (T\otimes I_2)\xtilde{}{t} $.
Rewriting~\eqref{eq:agent-dynamics-1} as a double integrator in state-space form 
with state $ \tilde{s}_j(\cdot) $ yields
\begin{equation}\label{eq:2n-order-system-state-space-1}
	\begin{array}{c}
		d\tilde{s}_j(t) = \left(F_0\tilde{s}_j(t)+F_{1j}\tilde{s}_j(t-1)\right)dt + Gd\noisebar{j}{t}, \\[5pt]
		 F_0 = \begin{bmatrix}
			0 & 1\\
			0 & -\gvel
		\end{bmatrix}, \ F_{1j} = \begin{bmatrix}
			0 & 0\\
			-\gvel\gpos_j & 0
		\end{bmatrix}, \ G = \begin{bmatrix}
			0\\
			1
		\end{bmatrix},
	\end{array}
\end{equation}
Stability of~\eqref{eq:multi-agent-state-space}
is equivalent to that of~\eqref{eq:2n-order-system-state-space-1} for all $ j $.
In the following, we drop the subscript $ j $ for the sake of readability.
%The first subsystem is stable for any $ \gvel > 0 $.
For positive eigenvalues $ \gpos $,~\eqref{eq:2n-order-system-state-space-1} is mean-square asymptotically stable
if $ \alpha_0 < 0 $ and unstable if $ \alpha_0 > 0 $~\cite{wangBoundedness}, where the \emph{spectral abscissa} is defined as
\begin{equation}\label{eq:2nd-order-system-stability-condition}
	\alpha_0 \doteq \sup\left\lbrace\Re(z) : z\in \mathbb{C}, \ h(z) = 0 \right\rbrace,
\end{equation}
and the \emph{characteristic polynomial} of~\eqref{eq:2n-order-system-state-space-1} is
\begin{equation}\label{eq:2n-order-system-chacteristic-polynomial}
	\begin{aligned}
		h(z) &\doteq \det\left(zI - F_0 - F_{1}\e^{-z}\right) = z^2 + \gvel z + \gvel\gpos\e^{-z}.  % = z^2 + \gvel z + \gvel\gpos\e^{-z}
%			 &= (z^2+\gvel z)\e^{z} + \gvel\gpos
	\end{aligned}
\end{equation}
A sufficient and necessary condition for all roots of $ h(z) $ to lie in the open left-hand half-plane is derived in~\cite{BAPTISTINI1997259}.
%and rewritten below for the sake of convenience.
\begin{thm}[\!\!\protect{\cite[Theorem 2.1]{BAPTISTINI1997259}}]\label{thm:stable-roots}
	Let the 2-vectors $v(b)=\left(p b, q-b^{2}\right), w(b)=$ $(\cos b, \sin b), b \geq 0,$ be given. If $r>0,$ a necessary and sufficient condition for all roots of the equation $h(z)=(z^2+pz+q)\e^{z} + r=0$ to have negative real part is that the orthogonality condition $v(b) \cdot w(b)=0,$ with $b \in \cup_{k=0}^{\infty}(2 k \pi,(2 k+1) \pi),$ implies $|v(b)|>r$.
\end{thm}
%In virtue of the above result, 
From~\cref{thm:stable-roots},~\eqref{eq:2n-order-system-state-space-1}
is asymptotically stable if the following implication holds for $b \in \cup_{k=0}^{\infty}(2 k \pi,(2 k+1) \pi)$,
\begin{equation}\label{eq:stability-condition-1}
	\gvel b\cos b - b^2\sin b = 0 \implies \gvel^2b^2+b^4>\gvel^2\gpos^2.
\end{equation}
In view of $ b \ge 0 $ and $ \sin b \ge 0 $,~\eqref{eq:stability-condition-1}
leads to~\eqref{eq:cont-time-double-int-stability-condition} after standard algebraic manipulations,
where we replace $ b $ with $ \beta = \min b \in (0,\nicefrac{\pi}{2}) $.
The inequality can be rewritten as
\begin{equation}\label{eq:stability-condition-lambda}
	\gpos < \dfrac{\beta}{\sin\beta} \doteq \phi(\gvel),
\end{equation}
where the definition of $ \phi(\cdot) $ follows from the implicit function theorem
applied to $ F(\gvel,\beta) \doteq \beta\tan\beta - \gvel $,
which states that $ F(\gvel,\beta) = 0 $ if and only if $ \beta = \varphi(\gvel) $ and
%$ \varphi'(\cdot), \varphi''(\cdot) $ can be explicitly computed, with
\begin{equation}\label{eq:varphi-of-eta-derivative}
	\varphi'(\gvel) = \dfrac{\cos^2\left(\varphi(\gvel)\right)}{\varphi(\gvel) + \sin\left(\varphi(\gvel)\right)\cos\left(\varphi(\gvel)\right)}
\end{equation}
Tedious but straightforward calculations on the first and second derivatives %$ \phi'(\gvel) $ and $ \phi''(\gvel) $
show that $ \phi(\gvel) $ is concave increasing for any $ \gvel > 0 $.
The limits at $ 0 $ and $ +\infty $ can be easily computed
by noting that
\begin{equation}\label{eq:stability-condition-limits}
	\beta_0\doteq\varphi(0) = 0, \quad \beta_\infty\doteq\lim_{\gvel\rightarrow+\infty}\varphi(\gvel) = \dfrac{\pi}{2}.
\end{equation}

	%!TEX ROOT = ../centralized_vs_distributed.tex

\subsection{\titlecap{derivation of first-order reduced model for continuous-time double integrators}}\label{app:time-scale-separation}

We now show that subsystem~\eqref{eq:agent-dynamics-1}
can be approximated to first-order dynamics
%when the control input is sufficiently powerful.
when the gain $ \gvel $ is sufficiently high.
%We first rewrite~\eqref{eq:agent-dynamics-1} in state-space form:
%\begin{equation}\label{eq:2n-order-system-state-space}
%\begin{aligned}
%d\x{}{t} &= \z{}{t}dt\\
%d\z{}{t} &= \left(-\gvel \z{}{t} - \gvel\gpos \x{}{t-1}\right)dt + d\noise{}{t}
%\end{aligned}
%\end{equation}
Let us consider~\eqref{eq:2n-order-system-state-space-1} with state $ \tilde{s}(t) = [\xtilde{}{t},\ztilde{}{t}]^\top $.
Assume that the feedback gain $ \gvel $ is large,
so that the variable $ \ztilde{}{t} $ evolves faster than $ \xtilde{}{t} $.
We can then approximate the dynamics of $ \ztilde{}{t} $ 
by letting $ \xtilde{}{t-1} \equiv x_0 $ be constant overtime,
\begin{equation}\label{eq:z-dynamics-with-constant-x}
d\ztilde{}{t} = \left(-\gvel \ztilde{}{t} - \gvel\gpos x_0\right)dt + d\noise{}{t}.
\end{equation}
\cref{eq:z-dynamics-with-constant-x} defines a standard Ornstein–Uhlenbeck process,
\begin{equation}\label{eq:z-solution-constant-x}
\ztilde{}{t} \sim \gauss\left( \e^{-\gvel t}(\ztilde{}{0} + \gpos x_0) - \gpos x_0, \dfrac{1}{2\gvel}\left(1-\e^{-2\gvel t}\right) \right).
\end{equation}
In view of the time-scale separation,
we assume that~\eqref{eq:z-solution-constant-x} holds (with $ \xtilde{}{t-1} $ constant) till $ \ztilde{}{t} $ settles at steady state, %the limit, tends to
\begin{equation}\label{eq:z-solution-constant-x-limit-gvel}
\lim_{t \rightarrow +\infty}\ztilde{}{t} = \tilde{z}_{\infty} \sim \gauss\left( - \gpos x_0, \dfrac{1}{2\gvel} \right).
\end{equation}
Using~\eqref{eq:z-solution-constant-x-limit-gvel},
we now approximate the dynamics of $ \xtilde{}{t} $ %in~\eqref{eq:agent-dynamics-1}
as if $ \ztilde{}{t} $ reached the steady state instantaneously,
\begin{equation}\label{eq:x-dynamics-1st-order}
d\xtilde{}{t} \approx \tilde{z}_{\infty}dt = -\gpos\xtilde{}{t-1}dt + dn(t),
\end{equation}
where 
%the drift only contains the dominant term $ -\gpos\x{}{t-1} $ and
the diffusion is embedded into the Brownian noise $ n(t) $ with variance proportional to $ \nicefrac{1}{\gvel} $.
In particular, as $ \gvel \rightarrow +\infty $,
$ \tilde{z}_{\infty} \xrightarrow{a.s.} - \gpos x_0 $
and~\eqref{eq:x-dynamics-1st-order} tends to deterministic dynamics.
	%!TEX ROOT = ../centralized_vs_distributed.tex

\newpage
\subsection{\titlecap{Computation of suboptimal variance for continuous-time single integrators}}\label{app:cont-time-single-int-suboptimal-variance-computation}

%Consider the suboptimal gain $ \tilde{k}^* $ as per~\cref{prop:subopt-gain}.
The $ N $ suboptimal eigenvalues have expression (cf.~\cite{circulant})
\begin{equation}\label{eq:eigenvaluesCirculant}
	\tilde{\gpos}_j^* = 2\tilde{k}^* \left(n - \sum_{\ell=1}^n\cos\left(\dfrac{2\pi (j-1) \ell}{N}\right)\right), %, \ j=1,\dots,N
\end{equation}
which we write as $ \tilde{\gpos}_j^* = g_j(n)\tilde{k}^* $.
Being $ \tilde{k}^* = \tilde{\alpha}^*(n) \opteig $
according to~\cref{prop:subopt-gain},
we write $ \tilde{\gpos}_j^* = \tilde{c}_j^*(n)\opteig $ with $  \tilde{c}_j^*(n) \doteq g_j(n)\tilde{\alpha}^*(n) $.
Then, each subsystem~\eqref{eq:cont-time-single-int-subsystem} has variance %$ \sigma_{\textit{ss}}^{2,I}(\tilde{\gpos}_j^*) $ is
\begin{align}\label{eq:optimalVarianceLambda_i}
	\varx{\tilde{\lambda}_j^*}{I} = \dfrac{1+\sin(\tilde{\lambda}_j^* \taun)}{2\tilde{\lambda}_j^*\cos(\tilde{\lambda}_j^* \taun)}
	\overset{(i)}{=} \underbrace{\dfrac{1+\sin(\tilde{c}_j^*(n)\beta^*)}{2\tilde{c}_j^*(n)\beta^*\cos(\tilde{c}_j^*(n)\beta^*)}}_{\doteq\tilde{C}_{j}^*(n)}\taun 
%	= \tilde{C}_{j}^*(n)\taun,
%	\begin{split}
%		\varx{\tilde{\lambda}_j^*}{I} &= \dfrac{1+\sin(\tilde{\lambda}_j^* \taun)}{2\tilde{\lambda}_j^*\cos(\tilde{\lambda}_j^* \taun)} \\
%		%&= \dfrac{1+\sin(g(n)\tilde{c}^*(n)\opteig \taun)}{2g(n)\tilde{c}^*(n)\opteig\cos(g_n)\tilde{c}^*(n)\opteig \taun)} \\
%		&\overset{(i)}{=} \dfrac{1+\sin(\tilde{c}_j^*(n)\beta^*)}{2\tilde{c}_j^*(n)\beta^*\cos(\tilde{c}_j^*(n)\beta^*)}\taun = \tilde{C}_{j}^*(n)\taun,
%	\end{split}
\end{align}
where~\eqref{eq:optimal-variance-closed-form} is used in \textit{(i)}.
%$ a_i^*(n) \doteq g_i(n)\tilde{c}^*(n) $ and
%$ \tilde{C}_i^*(n) $ multiplies $ \taun $ in the third line of~\eqref{eq:optimalVarianceLambda_i}.
%The scalar variance can then be written as
%\begin{equation}\label{eq:scalarVarianceSumExplicit}
%\optvarx = \sum_{i=2}^N \tilde{C}_i^*(n)\taun = \tilde{C}^*(n)f(n)\tau_{\textit{min}}
%\end{equation}
%where $ \displaystyle \tilde{C}^*(n) \doteq \sum_{i=2}^N \tilde{C}_i^*(n) $.
	%!TEX ROOT = ../centralized_vs_distributed.tex

\subsection{\titlecap{stability conditions for discrete-time systems}}\label{app:disc-time-single-int-stability}

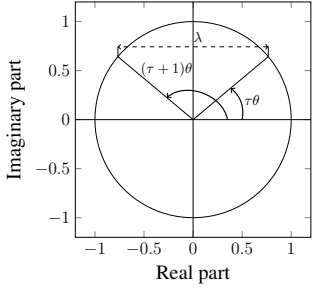
\begin{figure}
	\centering
	%!TEX ROOT = ../centralized_vs_distributed.tex

\begin{tikzpicture}[scale=.55]
	\begin{axis}[
		xlabel = Real part,
		ylabel = Imaginary part,
		axis equal image,
		ymin = -1.2,
		ymax = 1.2,
		xmin = -1.2,
		xmax = 1.2]
		\addplot [domain=-180:180, samples=100] ({cos(x)},{sin(x)});
		\draw (-1.2,0) -- (1.2,0);
		\draw (0,-1.2) -- (0,1.2);
		\draw (0,0) -- ({cos(40)},{sin(40)});
		\draw (0,0) -- ({-cos(40)},{sin(40)});
		\path[->,thick] (.5,0) edge[bend right] node [left] {}  ({.5*cos(40)},{.5*sin(40)});
		\path[->,thick] (.35,0) edge[bend right=60] node [left] {}  ({-.35*cos(40)},{.35*sin(40)}) ;
		\node at (.6,.2) {$ \delayn \theta $};
		\node at (-.27,.51) {$ (\delayn+1) \theta $};	
		\draw[<->,dashed] ({cos(40)},{sin(40)+.1}) -- ({-cos(40)},{sin(40)+.1});
		\draw ({cos(40)},{sin(40)+.1}) -- ({cos(40)},{sin(40)});
		\draw (-{cos(40)},{sin(40)+.1}) -- (-{cos(40)},{sin(40)});
		\node at (0.05,{sin(40)+.2}) {$ \gpos $};
	\end{axis}
\end{tikzpicture}
	\caption{A solution of~\eqref{eq:disc-time-single-int-uint-poles-system} in the complex plane.}
	\label{fig:disc-time-single-int-unit-poles}
\end{figure}

\myParagraph{\titlecap{General case}}
In the following, we replace $ \taun $ with $ \delayn $ % drop the subscript from $ \taun $
for the sake of readability.
For the single-integrator case, decoupling the error dynamics yields scalar subsystems of the form
\begin{equation}\label{eq:disc-time-single-int-decoupled}
	\xtilde{}{k+1} = \xtilde{}{k} - \gpos\xtilde{}{k-\tau} + \noisetilde{}{k}.
\end{equation}
The characteristic polynomial $ h(z) $ of~\eqref{eq:disc-time-single-int-decoupled} is obtained by
applying the lag operator $ z $
such that $ \xtilde{}{k}h(z) = \noisetilde{}{k} $,
\begin{equation}\label{eq:disc-time-single-int-characteristic-polinomial}
	h(z) = z - 1 + \gpos z^{-\tau}.
\end{equation}
Similarly, the double-integrator decoupled subsystems are
\begin{equation}\label{eq:disc-time-double-int-decoupled}
%	\xtilde{i}{k+1} = (2-\gvel)\xtilde{i}{k} - (1-\gvel)\xtilde{i}{k-1} - \gvel\gpos\xtilde{i}{k-\taun-1} + \noisetilde{i}{k}
	\begin{aligned}
		\xtilde{}{k+1} &= \xtilde{}{k} + \ztilde{}{k}\\
		\ztilde{}{k+1} &= (1-\gvel)\ztilde{}{k} - \gvel\gpos\xtilde{}{k-\tau} + \noisetilde{}{k},
	\end{aligned}
\end{equation}
with characteristic polynomial
\begin{equation}\label{eq:disc-time-double-int-characteristic-polinomial}
	h(z) = z-2+\gvel+(1-\gvel)z^{-1}+\gvel\gpos z^{-\tau-1}.
\end{equation}
For positive $ \gpos $, stability of~\eqref{eq:disc-time-single-int-decoupled}--\eqref{eq:disc-time-double-int-decoupled}
can be assessed via the Jury stability criterion,
which provides necessary and sufficient conditions for
the roots of~\eqref{eq:disc-time-single-int-characteristic-polinomial} and~\eqref{eq:disc-time-double-int-characteristic-polinomial}
to lie inside the unit circle %in the complex plane
in the form of inequalities involving the coefficients of $ h(z) $.
Being the latter polynomial in $ \gvel $ and $ \gpos $,
the Jury criterion %applied to~\eqref{eq:disc-time-single-int-model}--\eqref{eq:disc-time-double-int-model}
yields $ \Theta(N\tau) $ polynomial inequalities in the feedback gains,
which can be computed through standard software tools.

\myParagraph{Proof of~\cref{prop:disc-time-single-int-stability}}
\cref{eq:disc-time-single-int-characteristic-polinomial} can be studied as a root locus
by varying the gain $ \gpos $.
In particular, $ \gpos = 0 $ yields
a multiple root at $ z_1^* = 0 $ and a simple root at $ z_2^* = 1 $.
Negative values of $ \gpos $ are discarded as they push the latter outside the unit circle.
As $ \gpos $ increases,
the branches leave the unit ball along their asymptotes.
%Notice that, in view of the structure of the root locus,
The admissible values for $ \gpos $ are upper bounded by a threshold gain $ \gpos_{\textit{th}} $
beyond which some roots leave the unit ball.
In particular, we are interested in the minimum gain for which at least one root lies exactly on the unit circle.
Thus, we are looking for roots of~\eqref{eq:disc-time-single-int-characteristic-polinomial}
of the form $ z = \e^{j\theta} $,
\begin{equation}\label{eq:disc-time-single-int-unit-poles-eq}
	\e^{j(\delayn+1)\theta} - \e^{j\delayn\theta} + \gpos = 0.
\end{equation}
\cref{eq:disc-time-single-int-unit-poles-eq} can be equivalently written as the system
\begin{equation}\label{eq:disc-time-single-int-uint-poles-system}
	\begin{cases}
		\cos((\delayn+1)\theta) - \cos(\delayn\theta) + \gpos = 0\\
		\sin((\delayn+1)\theta) = \sin(\delayn\theta).
	\end{cases}
\end{equation}
\autoref{fig:disc-time-single-int-unit-poles} depicts a solution of~\eqref{eq:disc-time-single-int-uint-poles-system} for $ \sin(\delayn\theta) > 0 $.
The case $ \sin(\delayn\theta)<0 $ is analogous
and is omitted.
Further, the solution $ (\delayn+1)\theta = \delayn\theta $
can be discarded because it implies $ \gpos = 0 $ and thus prevents asymptotic stability.
%On the other hand, 
%Therefore, we only focus on the case $ \sin(\delayn\theta)>0 $. % depicted in the left box in~\autoref{fig:disc-time-single-int-unit-poles}.
From basic trigonometric arguments (c.f.~\autoref{fig:disc-time-single-int-unit-poles}),
the second equation in~\eqref{eq:disc-time-single-int-uint-poles-system} implies
\begin{equation}\label{eq:disc-time-single-int-theta}
	\delayn\theta + \dfrac{\theta}{2} = \dfrac{\pi}{2} + 2k\pi \ \longrightarrow \ \theta = \dfrac{\pi+4k\pi}{2\delayn+1}, \quad 
\end{equation}
where we impose $ \theta \in [0,\pi] $
and thus $ k\in\{0,\dots,\floor{\nicefrac{\delayn}{2}}\} $.
This includes all possible cases,
because the roots of~\eqref{eq:disc-time-single-int-characteristic-polinomial} come in complex conjugates pairs.
From~\eqref{eq:disc-time-single-int-theta},
the first equation in~\eqref{eq:disc-time-single-int-uint-poles-system},
and the fact $ \cos((\delayn+1)\theta) = - \cos(\delayn\theta) $,
we retrieve
\begin{equation}\label{eq:disc-time-single-int-unit-roots-gain}
	\gpos = 2\cos\left(\dfrac{\pi\delayn+4k\pi\delayn}{2\delayn+1}\right).
\end{equation}
The right-hand term in~\eqref{eq:disc-time-single-int-unit-roots-gain} is monotone increasing in $ k $.
Indeed, taking the argument of the cosine modulus $ 2\pi $ yields
\begin{equation}\label{eq:disc-time-single-int-angle-modulus}
	\dfrac{\pi\delayn+4k\pi\delayn}{2\delayn+1} \ \mod \ 2\pi = \dfrac{\pi\delayn-2k\pi}{2\delayn+1} \in \left[0,\dfrac{\pi}{2}\right),
\end{equation}
which is nonnegative and monotone decreasing in $ k $ for any $ \tau $.
Finally, the upper bound for the gain $ \gpos $ is given by
\begin{equation}\label{eq:disc-time-single-int-threshold-gain}
	\gpos_{\textit{th}} = \min_k2\cos\left(\dfrac{\pi\delayn+4k\pi\delayn}{2\delayn+1}\right) = 2\cos\left(\dfrac{\pi\delayn}{2\delayn+1}\right).
% 2\sin\left(\dfrac{\pi}{2}\dfrac{1}{2\delayn+1}\right)
\end{equation}

	\if0\mode
	%!TEX ROOT = ../centralized_vs_distributed.tex

\subsection{\titlecap{variance computation for discrete-time systems}}\label{app:disc-time-single-int-variance-explicit}

\myParagraph{Wiener--Kintchine Formula}
Given any fixed values of delay and feedback gains,
the steady-state variance $ \varx{\gpos}{I} $ or $ \varx{\gvel,\gpos}{II} $ of the decoupled subsystems can be computed numerically by
\begin{equation}\label{eq:disc-time-variance-integral}
	\dfrac{1}{2\pi}\int_{-\pi}^{+\pi}\dfrac{d\theta}{|h(\e^{j\theta})|^2},
\end{equation}
where the characteristic polynomial $ h(z) $
is~\eqref{eq:disc-time-single-int-characteristic-polinomial} or~\eqref{eq:disc-time-double-int-characteristic-polinomial}.

\myParagraph{Single Integrator Model}
The moment-matching method applied to subsystem~\eqref{eq:disc-time-single-int-decoupled} yields 
a linear system of equations in the variables $ (\rho_0,...,\rho_\delayn) $,
where $ \rho_t \doteq \mathbb{E}[\xtilde{}{k}\xtilde{}{k\pm t}] $:
\begin{subequations}\label{eq:disc-time-single-int-moment-matching-eqs}
	\begin{align}
		\rho_0 &= \mathbb{E}[\xtilde{}{k+1}^2] = \rho_0 + \gpos^2\rho_0 + 1 - 2\gpos\rho_\delayn \label{eq:disc-time-single-int-first-moment-eq}\\
		\rho_1 &= \mathbb{E}[\xtilde{}{k+1}\xtilde{}{k}] = \rho_0 - \gpos\rho_\delayn \label{eq:disc-time-single-int-yule-walker-eqs-1}\\
		&\hspace{2mm}\vdots \nonumber \\
		\rho_\delayn &= \rho_{\delayn-1} - \gpos\rho_1, \label{eq:disc-time-single-int-yule-walker-eqs-2}
	\end{align}
\end{subequations}
where~\eqref{eq:disc-time-single-int-yule-walker-eqs-1}--\eqref{eq:disc-time-single-int-yule-walker-eqs-2} are the Yule-Walker equations.
%associated to the decoupled subsystem~\eqref{eq:disc-time-single-int-decoupled}.
System~\eqref{eq:disc-time-single-int-moment-matching-eqs}
can be written compactly as $ A^{(\tau)}\rho = e_1 $, 
where $ \rho^\top = [\rho_0,\dots,\rho_{\delayn}]$,
$ e_1 $ is the canonical vector in $ \Real{\delayn+1} $ with nonzero first coordinate,
and $ A^{(\tau)}\in\Real{(\delayn+1)\times(\delayn+1)} $ gathers all coefficients of equations in~\eqref{eq:disc-time-single-int-moment-matching-eqs}.
%with
%\begin{equation}\label{eq:explicit-variance-matrix-A}
%A^{(\tau)} = \begin{bmatrix}
%-\gpos^2 &   		&     		& 		 &   	  & 2\gpos\\
%1 		 & -1 		&     		& 		 &    	  & -\gpos\\
%& 	\ddots	& \ddots	&  		 & \iddots&  \\
%& 			& 			&		 & 		  &  \\
%& 			& 			&		 & 		  &  \\
%& -\gpos	& 			& 		 & 1 	  & -1
%\end{bmatrix}.
%\end{equation}
%In particular, when $ \delayn $ is odd, the $ (\ceil{\nicefrac{\delayn}{2}} + 1) $-th row is
%\begin{equation}\label{eq:app-explicit-variance-matrix-A-mid-row-odd}
%	\left[\!\begin{array}{cccccccc}
%	0 & \dots & 0 & 1 & -1-\gpos & 0 & \dots & 0
%	\end{array}\!\right],
%\end{equation}
%while, when $ \delayn $ is even, the $ (\nicefrac{\delayn}{2} + 2) $-th row is
%\begin{equation}\label{eq:app-explicit-variance-matrix-A-mid-row-even}
%	\left[\!\begin{array}{cccccccc}
%		0 & \dots & 0 & 1-\gpos & -1 & 0 & \dots & 0
%	\end{array}\!\right].
%\end{equation}
%where we highlight the dependence on the delay in view of the recursive characterization of $ \rho_0 $.
It can be seen that $ A^{(\tau)} $ is full rank for all $ \delayn \ge 1 $ and thus~\eqref{eq:disc-time-single-int-moment-matching-eqs} has a unique solution.
In particular, we are interested in the autocorrelation $ \rho_0 = \varx{\gpos}{I} $,
which is given by
the ratio between the minor associated with the top-left element of $ A^{(\tau)} $,
named $ n_\delayn \doteq M^{(\tau)}_{1,1} $, and the determinant $ d_\delayn \doteq \det(A^{(\tau)}) $.
Specifically, $ \rho_0 $ is a rational function in $ \gpos $
and can be computed in closed form by a symbolic solver
given any value of $ \delayn $.

Further, $ n_{\delayn} $ and $ d_{\delayn} $ can be explicitly computed %recursively via an inductive argument on the delay $ \delayn $, 
by leveraging a recursive nested structure of the matrix $ A $.
%\begin{equation}\label{eq:recursive-matrix}
%	A^{(\delayn)} = \tikz[baseline=(M.west)]{%
%		\node[matrix of math nodes,matrix anchor=west,left delimiter={[},right delimiter={]},ampersand replacement=\&] (M) {%
%		-\gpos^2\&			\&			\&							\& 							\&						 \&	-2\gpos		\\
%			1 	\& -1  		\&    		\&							\& 							\& 						 \&	  			\\
%				\&  1  		\&  -1		\&	  						\&      					\& 						 \& -\gpos   	\\
%				\&    		\&   1 		\& \textcolor{white}{t1}	\& 							\& 						 \&				\\
%				\&			\&			\&							\& \tilde{A}^{(\delayn-4)}	\&						 \&				\\
%				\& 			\&	  		\& 							\&	 						\& \textcolor{white}{t1} \&				\\
%				\& 		    \&   -\gpos	\&							\&							\& 1					 \&	-1			\\
%				\& 	-\gpos	\&	  		\&							\&							\&						 \&	1			\\			
%		};
%		\node[draw,fit=(M-3-3)(M-7-7),inner sep=-1pt] {};
%		\node[draw,fit=(M-4-4)(M-6-6),inner sep=-1pt] {};
%	}.
%\end{equation}
%where $ \tilde{A}^{(\delayn)} $ is the submatrix of $ A^{(\delayn)} $ obtained by removing its first row and column
%such that $ M^{(\delayn)}_{1,1} = \det(\tilde{A}^{(\delayn)}) $,
%and the matrices $ \tilde{A}^{(\delayn-2)} $ and $ \tilde{A}^{(\delayn-4)} $ are framed in~\eqref{eq:recursive-matrix}.
The solution obeys the following recursive expression in $ \delayn $:
%\begin{prop}\label{prop:disc-time-single-int-variance-explicit}
\begin{subequations}\label{eq:disc-time-single-int-variance-explicit}
	\begin{gather}
		n_\delayn = \begin{dcases}
			(-1-\gpos)n_{\delayn-1} + \tilde{n}_{\delayn-1} & \mbox{if } \delayn \mbox{ odd}\\
			-(1-\gpos)n_{\delayn-1} - \gpos\tilde{n}_{\delayn-1} & \mbox{if }\delayn \mbox{ even},\\
		\end{dcases} \label{eq:disc-time-single-int-variance-explicit-numerator}\\
		\tilde{n}_{\delayn} = (2-\gpos^2)\tilde{n}_{\delayn-2} - \tilde{n}_{\delayn-4}, \label{eq:disc-time-single-int-variance-explicit-numerator-rem}\\
		d_\delayn = d_{\delayn-2} - \gpos^2\left(n_\delayn+n_{\delayn-2}\right), \label{eq:disc-time-single-int-variance-explicit-denominator}\\
		\tilde{n}_{-3} = -1+\gpos^2, \ \tilde{n}_{-2} = \gpos^2, \ \tilde{n}_{-1} = -1, \ \tilde{n}_0 = 0, \\
		n_{-1} = 0, \ n_0 = 1, \ d_{-1} = -2\gpos, \ d_0 = 2\gpos-\gpos^2.
	\end{gather}
\end{subequations}
Detailed derivation of~\eqref{eq:disc-time-single-int-variance-explicit}
is given in the technical report~\cite{2021arXiv210900359B}.
Given $ \delayn $,
convexity of $ \rho_0 $ in $ \gpos $ can be assessed by checking the sign
of the second derivative in the stability region.
This reduces to a system of inequalities
which can be solved, \eg by \texttt{solve\_rational\_inequalities} in Python.
The variance was proved strictly convex for all tried delays.

\myParagraph{\titlecap{double integrator model}}
System~\eqref{eq:disc-time-double-int-decoupled} yields the following
$ \delayn+2 $ coupled moment-matching equations,
%The moment-matching system associated with~\eqref{eq:disc-time-double-int-decoupled}
%has $ \delayn+2 $ variables $ (\rho_0,\dots,\rho_{\delayn+1}) $ and is composed of the following equations:
\begin{subequations}\label{eq:disc-time-double-int-moment-matching-eqs}
	\begin{align}
		\begin{split}
			\rho_0 &= [(2-\gvel)^2 + (1-\gvel)^2 + \gvel^2\gpos^2]\rho_0 - 2(2-\gvel)(1-\gvel)\rho_1 \\
			& + 2(1-\gvel)\gvel\gpos\rho_\delayn - 2(2-\gvel)\gvel\gpos\rho_{\delayn+1} + 1 \label{eq:disc-time-double-int-first-moment-eq}
		\end{split}\\
		\rho_1 &= (2-\gvel)\rho_0 - (1-\gvel)\rho_1 - \gvel\gpos\rho_{\delayn+1} \label{eq:disc-time-double-int-yule-walker-eqs-1}\\
		&\hspace{2mm}\vdots \nonumber \\
		\rho_{\delayn+1} &= (2-\gvel)\rho_\delayn - (1-\gvel)\rho_{\delayn-1} - \gvel\gpos\rho_1, \label{eq:disc-time-double-int-yule-walker-eqs-2}
	\end{align}
\end{subequations}
with~\eqref{eq:disc-time-double-int-yule-walker-eqs-1}--\eqref{eq:disc-time-double-int-yule-walker-eqs-2} the associated Yule-Walker equations.
%associated with~\eqref{eq:disc-time-double-int-decoupled}.
Analogous analysis to single-integrator model can be performed.
	\else
	%!TEX ROOT = ../centralized_vs_distributed.tex

\subsection{\titlecap{variance computation for discrete-time systems}}\label{app:disc-time-single-int-variance-explicit}

\myParagraph{Wiener--Kintchine Formula}
Given any fixed values of delay and feedback gains,
the steady-state variance $ \varx{\gpos}{I} $ or $ \varx{\gvel,\gpos}{II} $ of the decoupled subsystems can be computed numerically by
\begin{equation}\label{eq:disc-time-variance-integral}
	\dfrac{1}{2\pi}\int_{-\pi}^{+\pi}\dfrac{d\theta}{|h(\e^{j\theta})|^2},
\end{equation}
where the characteristic polynomial $ h(z) $
is~\eqref{eq:disc-time-single-int-characteristic-polinomial} or~\eqref{eq:disc-time-double-int-characteristic-polinomial}.

\myParagraph{\titlecap{single integrator model}}
The moment-matching method applied to the subsystem~\eqref{eq:disc-time-single-int-decoupled} yields 
a linear system of equations in the variables $ (\rho_0,...,\rho_\delayn) $,
where $ \rho_t \doteq \mathbb{E}[\xtilde{}{k}\xtilde{}{k\pm t}] $:
\begin{subequations}\label{eq:disc-time-single-int-moment-matching-eqs}
	\begin{align}
		\rho_0 &= \mathbb{E}[\xtilde{}{k+1}^2] = \rho_0 + \gpos^2\rho_0 + 1 - 2\gpos\rho_\delayn \label{eq:disc-time-single-int-first-moment-eq}\\
		\rho_1 &= \mathbb{E}[\xtilde{}{k+1}\xtilde{}{k}] = \rho_0 - \gpos\rho_\delayn \label{eq:disc-time-single-int-yule-walker-eqs-1}\\
		&\hspace{2mm}\vdots \nonumber \\
		\rho_\delayn &= \rho_{\delayn-1} - \gpos\rho_1, \label{eq:disc-time-single-int-yule-walker-eqs-2}
	\end{align}
\end{subequations}
where~\eqref{eq:disc-time-single-int-yule-walker-eqs-1}--\eqref{eq:disc-time-single-int-yule-walker-eqs-2} are the Yule-Walker equations.
%associated to the decoupled subsystem~\eqref{eq:disc-time-single-int-decoupled}.
System~\eqref{eq:disc-time-single-int-moment-matching-eqs}
can be written compactly as $ A^{(\tau)}\rho = e_1 $, where
$ \rho^\top = [\rho_0,\dots,\rho_{\delayn}]$,
$ e_1 $ is the canonical vector in $ \Real{\delayn+1} $ with nonzero first coordinate and $ A^{(\tau)}\in\Real{(\delayn+1)\times(\delayn+1)} $ with
\begin{equation}\label{eq:explicit-variance-matrix-A}
	A^{(\tau)} = \begin{bmatrix}
		-\gpos^2 &   		&     		& 		 &   	  & 2\gpos\\
		1 		 & -1 		&     		& 		 &    	  & -\gpos\\
		& 	\ddots	& \ddots	&  		 & \iddots&  \\
		& 			& 			&		 & 		  &  \\
		& 			& 			&		 & 		  &  \\
		& -\gpos	& 			& 		 & 1 	  & -1
	\end{bmatrix}.
\end{equation}
In particular, when $ \delayn $ is odd, the $ (\ceil{\nicefrac{\delayn}{2}} + 1) $-th row is
\begin{equation}\label{eq:app-explicit-variance-matrix-A-mid-row-odd}
	\left[\!\begin{array}{cccccccc}
		0 & \dots & 0 & 1 & -1-\gpos & 0 & \dots & 0
	\end{array}\!\right],
\end{equation}
while, when $ \delayn $ is even, the $ (\nicefrac{\delayn}{2} + 2) $-th row is
\begin{equation}\label{eq:app-explicit-variance-matrix-A-mid-row-even}
	\left[\!\begin{array}{cccccccc}
		0 & \dots & 0 & 1-\gpos & -1 & 0 & \dots & 0
	\end{array}\!\right].
\end{equation}
%where we highlight the dependence on the delay in view of the recursive characterization of $ \rho_0 $.
Notice that $ A^{(\tau)} $ is full rank for all $ \delayn \ge 1 $ and thus~\eqref{eq:disc-time-single-int-moment-matching-eqs} can be solved uniquely.
In particular, we are interested in the autocorrelation $ \rho_0 = \varx{\gpos}{I} $,
which is given by
the ratio between the minor associated with the top-left element of $ A^{(\tau)} $,
named $ n_\delayn \doteq M^{(\tau)}_{1,1} $, and the determinant $ d_\delayn \doteq \det(A^{(\tau)}) $.
Specifically, $ \rho_0 $ is a rational function in $ \gpos $
and can be computed in closed form by a symbolic solver
given any value of $ \delayn $.

Further, $ n_{\delayn} $ and $ d_{\delayn} $ can be computed %recursively via an inductive argument on the delay $ \delayn $, 
by leveraging the following nested structure of the matrix $ A^{(\delayn)} $:
\begin{equation}\label{eq:recursive-matrix}
	A^{(\delayn)} = \tikz[baseline=(M.west)]{%
		\node[matrix of math nodes,matrix anchor=west,left delimiter={[},right delimiter={]},ampersand replacement=\&] (M) {%
			-\gpos^2\&			\&			\&							\& 							\&						 \&	-2\gpos		\\
			1 	\& -1  		\&    		\&							\& 							\& 						 \&	  			\\
			\&  1  		\&  -1		\&	  						\&      					\& 						 \& -\gpos   	\\
			\&    		\&   1 		\& \textcolor{white}{t1}	\& 							\& 						 \&				\\
			\&			\&			\&							\& \tilde{A}^{(\delayn-4)}	\&						 \&				\\
			\& 			\&	  		\& 							\&	 						\& \textcolor{white}{t1} \&				\\
			\& 		    \&   -\gpos	\&							\&							\& 1					 \&	-1			\\
			\& 	-\gpos	\&	  		\&							\&							\&						 \&	1			\\			
		};
		\node[draw,fit=(M-3-3)(M-7-7),inner sep=-1pt] {};
		\node[draw,fit=(M-4-4)(M-6-6),inner sep=-1pt] {};
	},
\end{equation}
where $ \tilde{A}^{(\delayn)} $ is the submatrix of $ A^{(\delayn)} $ obtained by removing its first row and column
such that $ M^{(\delayn)}_{1,1} = \det(\tilde{A}^{(\delayn)}) $,
and the matrices $ \tilde{A}^{(\delayn-2)} $ and $ \tilde{A}^{(\delayn-4)} $ are framed in~\eqref{eq:recursive-matrix}.

The solution obeys the following recursive expression in $ \delayn $:
%\begin{prop}\label{prop:disc-time-single-int-variance-explicit}
\begin{subequations}\label{eq:disc-time-single-int-variance-explicit}
	\begin{gather}
		n_\delayn = \begin{dcases}
			(-1-\gpos)n_{\delayn-1} + \tilde{n}_{\delayn-1} & \mbox{if } \delayn \mbox{ odd}\\
			-(1-\gpos)n_{\delayn-1} - \gpos\tilde{n}_{\delayn-1} & \mbox{if }\delayn \mbox{ even},\\
		\end{dcases} \label{eq:disc-time-single-int-variance-explicit-numerator}\\
		\tilde{n}_{\delayn} = (2-\gpos^2)\tilde{n}_{\delayn-2} - \tilde{n}_{\delayn-4}, \label{eq:disc-time-single-int-variance-explicit-numerator-rem}\\
		d_\delayn = d_{\delayn-2} - \gpos^2\left(n_\delayn+n_{\delayn-2}\right), \label{eq:disc-time-single-int-variance-explicit-denominator}\\
		\tilde{n}_{-3} = -1+\gpos^2, \ \tilde{n}_{-2} = \gpos^2, \ \tilde{n}_{-1} = -1, \ \tilde{n}_0 = 0, \\
		n_{-1} = 0, \ n_0 = 1, \ d_{-1} = -2\gpos, \ d_0 = 2\gpos-\gpos^2.
	\end{gather}
\end{subequations}

\cref{eq:disc-time-single-int-variance-explicit} can be proved by an inductive argument
on the delay $ \delayn $.

\myParagraph{Numerator}
We demonstrate the formula for odd delays $ \delayn = 2k+1, k\in\mathbb{N} $.
The other case can be obtained similarly and is thus omitted. % in the interest of space.\\

Let us consider the submatrix $ \tilde{A}^{(\delayn)} \in\Real{\delayn\times\delayn} $
obtained by removing the first row and column of $ A $,
such that $ n_\delayn = \det(\tilde{A}^{(\delayn)}) $.
Replacing the $ (\floor{\nicefrac{\delayn}{2}}) $-th column with the sum of
$ (\floor{\nicefrac{\delayn}{2}}) $-th and $ (\ceil{\nicefrac{\delayn}{2}}) $-th columns yields
\begin{equation}\label{eq:numerator-1}
	\det\left(\tilde{A}^{(\delayn)}\right) = \left|\begin{array}{c|c|c}
		\tilde{A}_{11}^{(\delayn-1)} &  & \tilde{A}_{12}^{(\delayn-1)} \\
		\hline
		\begin{array}{ccc} \dots & 0 & -\gpos \end{array} & -1-\gpos & \\
		\hline
		\tilde{A}_{21}^{(\delayn-1)} & \begin{array}{c} 1 \\ 0 \\ \vdots \end{array} & \tilde{A}_{22}^{(\delayn-1)}
	\end{array}\right|,
\end{equation}
from which it follows $ n_\delayn = (-1-\gpos)n_{\delayn-1} - \det(R^{(\delayn)}) $ where $ R^{(\delayn)}\in\Real{(\delayn-1)\times(\delayn-1)} $
and the base case is $ n_1 = -1-\gpos $.
This expression corresponds to~\eqref{eq:disc-time-single-int-variance-explicit-numerator} with $ \tilde{n}_{\delayn-1} = -\det(R^{(\delayn)}) $.
Manipulations of the second term yield a further recursive expression for $ \tilde{n}_{\delayn-1} $.
Let us write
\begin{equation}\label{eq:numerator-2}
	\det\left(R^{(\delayn)}\right) =	\tikz[baseline=(M.west)]{%
		\node[matrix of math nodes,matrix anchor=west,left delimiter=|,right delimiter=|,ampersand replacement=\&] (M) {%
			-1 		\&    		\& 	  \& 	 			\& 				\&   		\& -\gpos\\
			1 		\& -1  		\&    \& 				\&  			\& -\gpos 	\&  \\
			\&  1  		\&  \textcolor{white}{t1}  \&       			\&    			\& 			\& \\
			\&    		\&    \& R^{(\delayn-4)}	\&	  			\& 			\& \\
			\& 			\&	  \&	 			\&	\textcolor{white}{t2}			\&			\&	\\
			\& \gpos    \&    \&				\&		1		\&	-1		\&	\\
			-\gpos	\& 			\&	  \&				\&				\& 	1		\& -1\\			
		};
		\node[draw,fit=(M-2-2)(M-6-6),inner sep=-1pt] {};
		\node[draw,fit=(M-3-3)(M-5-5),inner sep=-1pt] {};
	},
\end{equation}
where the two inner boxes highlight $ R^{(\delayn-2)} $ and $ R^{(\delayn-4)} $, respectively.
Straightforward calculations yield
\begin{multline}\label{eq:numerator-3}
	\det\left(R^{(\delayn)}\right) = \det\left(R^{(\delayn-2)}\right) + \\
	\gpos\tikz[baseline=(M.west)]{%
		\node[matrix of math nodes,matrix anchor=west,left delimiter=|,right delimiter=|,ampersand replacement=\&] (M) {%
			1 			\& -1  		\&    						\& 					\&  						\& -\gpos 	 \\
			\&  1  		\&  \textcolor{white}{t1}  	\&       			\&    						\& 			 \\
			\&    		\&    						\& R^{(\delayn-4)}	\&	  						\& 			 \\
			\& 			\&	  						\&	 				\&	\textcolor{white}{t2}	\&			 \\
			\& -\gpos   \&    						\&					\&		1					\&	-1		 \\
			-\gpos		\& 			\&	  						\&					\&							\& 	1		 \\			
		};
		\node[draw,fit=(M-1-2)(M-5-6),inner sep=-1pt] {};
		\node[draw,fit=(M-2-3)(M-4-5),inner sep=-1pt] {};
	}.
\end{multline}
The determinant in the second addend is computed as
\begin{equation}\label{eq:numerator-4}
	-\gpos\det\left(R^{(\delayn-2)}\right) + 
	\tikz[baseline=(M.west)]{%
		\node[matrix of math nodes,matrix anchor=west,left delimiter=|,right delimiter=|,ampersand replacement=\&] (M) {%
			1  		 \&  \textcolor{white}{t1}  \&       			\&    						\\
			\&    						\& R^{(\delayn-4)}	\&	  						\\
			\&	  						\&	 				\&	\textcolor{white}{t2}	\\
			-\gpos   \&    						\&					\&		1					\\
		};
		\node[draw,fit=(M-1-2)(M-3-4),inner sep=-1pt] {};
	},
\end{equation}
and the second addend in the above equation has the same structure as
the determinant in the second addend in~\eqref{eq:numerator-3}.
Thus, an easy inductive argument proves
\begin{multline}\label{eq:numerator-rem-recursive}
	\det\left(R^{(\delayn)}\right) = \det\left(R^{(\delayn-2)}\right) + \gpos\left(-\gpos\det\left(R^{(\delayn-2)}\right)\right.\\
	\left.-\gpos\det\left(R^{(\delayn-4)}\right) - \dots - \gpos\det\left(R^{(3)}\right)- \gpos\right),
\end{multline}
where the base case is $ \det\left(R^{(3)}\right) = -\gpos^2 $.
\cref{eq:disc-time-single-int-variance-explicit-numerator-rem} is retrieved by noting
\begin{multline}\label{eq:numerator-rem-recursive-1}
	\det\left(R^{(\delayn-2)}\right) - \left(1-\gpos^2\right)\det\left(R^{(\delayn-4)}\right) = \\
	\gpos\left(-\gpos\det\left(R^{(\delayn-6)}\right) - \dots -\gpos\right),
\end{multline}
and thus the tail of the infinite summation in~\eqref{eq:numerator-rem-recursive} can be replaced
by the left-hand term in~\eqref{eq:numerator-rem-recursive-1}.

\myParagraph{Denominator}
The denominator of $ \rho_0 $ is computed as the determinant of $ A $.
Let $ A^{(\delayn)} \doteq A $, from~\eqref{eq:explicit-variance-matrix-A} we get
\begin{multline}\label{eq:denominator-1}
	\det\left(A^{(\delayn)}\right) = -\gpos^2M_{1,1}^{(\delayn)} -\\
	2\gpos\tikz[baseline=(M.west)]{%
		\node[matrix of math nodes,matrix anchor=west,left delimiter=|,right delimiter=|,ampersand replacement=\&] (M) {%
			1 	\& -1  		\&    		\&							\& 							\& 						 \&	  			\\
			\&  1  		\&  -1		\&	  						\&      					\& 						 \& -\gpos   	\\
			\&    		\&   1 		\& \textcolor{white}{t1}	\& 							\& 						 \&				\\
			\&			\&			\&							\& \tilde{A}^{(\delayn-4)}	\&						 \&				\\
			\& 			\&	  		\& 							\&	 						\& \textcolor{white}{t1} \&				\\
			\& 		    \&   -\gpos	\&							\&							\& 1					 \&	-1			\\
			\& 	-\gpos	\&	  		\&							\&							\&						 \&	1			\\			
		};
		\node[draw,fit=(M-2-3)(M-6-7),inner sep=-1pt] {};
		\node[draw,fit=(M-3-4)(M-5-6),inner sep=-1pt] {};
	},
\end{multline}
where $ \tilde{A}^{(\delayn-2)} $ and $ \tilde{A}^{(\delayn-4)} $ are framed in the second addend above.
The latter can be computed as the following sum,
\begin{equation}\label{eq:denominator-2}
	\gpos M_{1,1}^{(\delayn-2)} + \tikz[baseline=(M.west)]{%
		\node[matrix of math nodes,matrix anchor=west,left delimiter=|,right delimiter=|,ampersand replacement=\&] (M) {%
			1  		\&  -1				\&	  						\&      								\& 							\\
			\&   1 				\& \textcolor{white}{t1}	\& 										\&							\\
			\&					\&							\& \tilde{A}^{(\delayn-4)}				\&							\\
			\&	  				\& 							\&	 									\& 	\textcolor{white}{t1}	\\
			\&   -\gpos			\&							\&										\& 1						\\
		};
		\node[draw,fit=(M-2-3)(M-4-5),inner sep=-1pt] {};
	},
\end{equation}
where the same structure is repeated recursively in the second addend above.
Thus, an easy inductive argument proves
\begin{equation}\label{eq:denominator-recursirve}
	d_\delayn = -\gpos^2n_\delayn - 2\gpos\left(\gpos n_{\delayn-2} + \gpos n_{\delayn-4} + \dots + \gpos n_1 + 1\right),
\end{equation}
where the base case is $ d_1 = -\gpos^2(-1-\gpos)-2\gpos $.
\cref{eq:disc-time-single-int-variance-explicit-denominator} is retrieved by noting
\begin{multline}\label{eq:denominator-recursirve-1}
	-2\gpos\left(\gpos n_{\delayn-2} + \gpos n_{\delayn-4} + \dots + 1\right) = \\
	-\gpos^2 n_{\delayn-2} -\gpos^2 n_{\delayn-2} - 2\gpos\left(\gpos n_{\delayn-4} + \dots + 1\right) = \\
	-\gpos^2 n_{\delayn-2} + d_{\delayn-2}.
\end{multline}

Given $ \delayn $,
convexity of $ \rho_0 $ in $ \gpos $ can be assessed by checking the sign
of the second derivative in the stability region.
This reduces to a system of inequalities
which can be solved, \eg by \texttt{solve\_rational\_inequalities} in Python.
The variance was proved strictly convex for all tried delays.

\myParagraph{\titlecap{double integrator model}}
The moment-matching system associated with~\eqref{eq:disc-time-double-int-decoupled}
has $ \delayn+2 $ variables $ (\rho_0,\dots,\rho_{\delayn+1}) $ and is composed of the following equations:
\begin{subequations}\label{eq:disc-time-double-int-moment-matching-eqs}
	\begin{align}
		\begin{split}
			\rho_0 &= (2-\gvel)^2\rho_0 + (1-\gvel)^2\rho_0 + \gvel^2\gpos^2\rho_0 + 1 \\
			&- 2(2-\gvel)(1-\gvel)\rho_1 - 2(2-\gvel)\gvel\gpos\rho_{\delayn+1} \\
			&+ 2(1-\gvel)\gvel\gpos\rho_\delayn
			\label{eq:disc-time-double-int-first-moment-eq}
		\end{split}\\
		\rho_1 &= (2-\gvel)\rho_0 - (1-\gvel)\rho_1 - \gvel\gpos\rho_{\delayn+1} \label{eq:disc-time-double-int-yule-walker-eqs-1}\\
		\rho_2 &= (2-\gvel)\rho_1 - (1-\gvel)\rho_0 - \gvel\gpos\rho_{\delayn}\\
		&\hspace{2mm}\vdots \nonumber \\
		\rho_{\delayn+1} &= (2-\gvel)\rho_\delayn - (1-\gvel)\rho_{\delayn-1} - \gvel\gpos\rho_1, \label{eq:disc-time-double-int-yule-walker-eqs-2}
	\end{align}
\end{subequations}
where~\eqref{eq:disc-time-double-int-yule-walker-eqs-1}--\eqref{eq:disc-time-double-int-yule-walker-eqs-2} are the Yule-Walker equations
associated with~\eqref{eq:disc-time-double-int-decoupled}.
Analogous considerations to the single-integrator model can be done in this case.
	\fi
	
	%!TEX root = ../centralized_vs_distributed.tex

\begin{IEEEbiography}[{\includegraphics[width=1in,height=1.25in,clip,keepaspectratio]{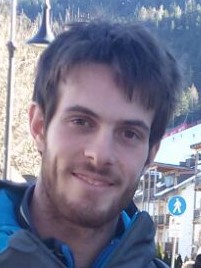}}]{Luca Ballotta}
	received his Master's Degree in Automation 
Engineering in 2019 from the University of Padova, where he is currently pursuing the Ph.D. degree in Information Engineering.
He was Visiting Student at the Massachusetts Institute of Technology in 2020 and 2022.
His research interests include networked control systems subject to resource constraints, 
resilient distributed control,
and learning-based safe control. 
He was awarded with the Young Author Prize at the 2020 IFAC World Congress.

\end{IEEEbiography}

\begin{IEEEbiography}[{\includegraphics[width=1in,height=1.25in,clip,keepaspectratio]{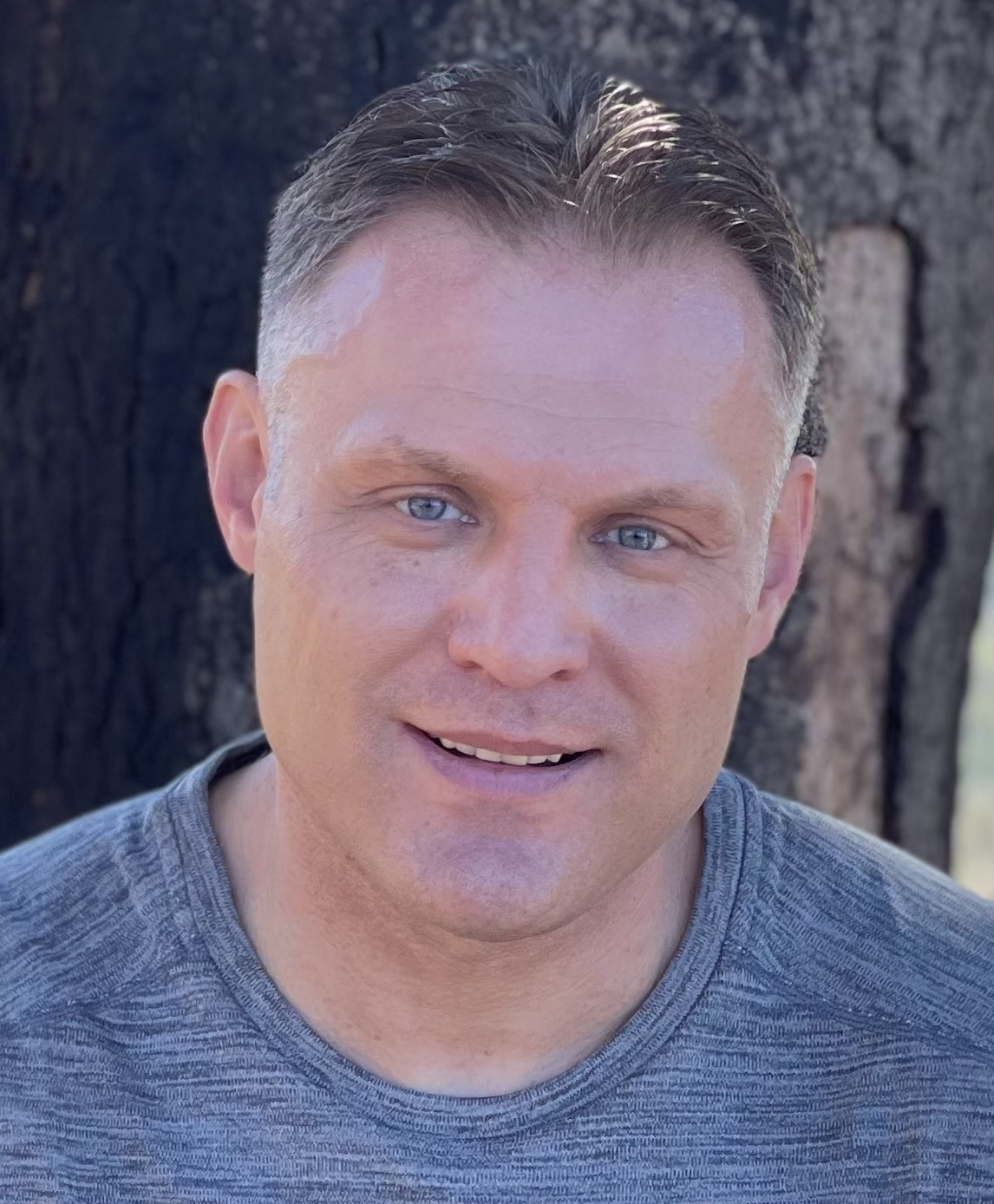}}]{Mihailo~R.~Jovanovi\'c}
	(Fellow IEEE) received the PhD degree in mechanical engineering from the University of California at Santa Barbara, Santa Barbara, CA, USA, in 2004. 
	
	He is currently a Professor in the Ming Hsieh Department of Electrical and Computer Engineering and the Founding Director of the Center for Systems and Control at the University of Southern California, Los Angeles, CA, USA. He was a faculty member in the Department of Electrical and Computer Engineering at the University of Minnesota, Twin Cities, MN, USA, from 2004 until 2017, and has held visiting positions with Stanford University, the Institute for Mathematics and its Applications, the Simons Institute for the Theory of Computing, and the University of Belgrade. 
	
	Prof.\ Jovanovi\'c received a CAREER Award from the National Science Foundation in 2007, the George S. Axelby Outstanding Paper Award from the IEEE Control Systems Society in 2013, and the Distinguished Alumnus Award from the Department of Mechanical Engineering, University of California at Santa Barbara, in 2014. Papers of his students were finalists for the Best Student Paper Award at the American Control Conference in 2007 and 2014. He is a Fellow of the American Physical Society. 

\end{IEEEbiography}

\begin{IEEEbiography}[{\includegraphics[width=1in,height=1.25in,clip,keepaspectratio]{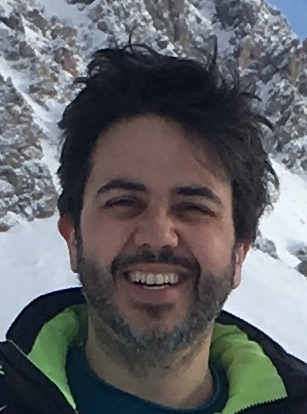}}]{Luca Schenato}
	(Fellow IEEE) received the Dr. Eng. degree in electrical engineering from the University of Padova in 1999 and the Ph.D. degree in Electrical Engineering and Computer Sciences from the UC Berkeley, in 2003. He held a post-doctoral position in 2004 and a visiting professor position in 2013-2014 at U.C. Berkeley. Currently he is Full Professor at the Information Engineering Department at the University of Padova. His interests include networked control systems, multi-agent systems, wireless sensor networks, distributed optimisation and synthetic biology. Luca Schenato has been awarded the 2004 Researchers Mobility Fellowship by the Italian Ministry of Education, University and Research (MIUR), the 2006 Eli Jury Award in U.C. Berkeley and the EUCA European Control Award in 2014, and IEEE Fellow in 2017. He served as Associate Editor for IEEE Trans. on Automatic Control from 2010 to 2014 and he is he is currently Senior Editor for IEEE Trans. on Control of Network Systems and Associate Editor for Automatica.

\end{IEEEbiography}
	
\end{document}